\documentclass[mathpazo,online]{jml}
\usepackage[toc,page]{appendix}

\usepackage{bbm}
\usepackage{multirow,multicol,tabularx}
\usepackage{tikz,pgfplots}
\usetikzlibrary{patterns}
\pgfplotsset{compat=newest}
\usepackage{hyperref}
\usepackage{float,cleveref}

\usepackage{caption}
\usepackage{graphicx}

\newcommand{\Psid}{\Psi^\dagger}
\newcommand{\cU}{\mathcal{U}}
\newcommand{\cV}{\mathcal{V}}
\newcommand{\cH}{\mathcal{H}}

\newcommand{\R}{\mathbb{R}}
\newcommand{\C}{\mathbb{C}}
\newcommand{\Z}{\mathbb{Z}}
\newcommand{\N}{\mathbb{N}}
\newcommand{\E}{\mathbb{E}}
\newcommand{\1}{\mathbbm{1}}
\newcommand{\cR}{\mathcal{R}}
\newcommand{\cQ}{\mathcal{Q}}
\newcommand{\cL}{\mathcal{L}}
\newcommand{\cK}{\mathcal{K}}
\newcommand{\cF}{\mathcal{F}}

\newcommand{\bigO}{\mathcal{O}}

\definecolor{darkred}{rgb}{.7,0,0}
\definecolor{darkgreen}{rgb}{0,.7,0}
\definecolor{darkblue}{rgb}{0.13,0.44,.72}
\newcommand{\eq}[1]{{\color{black}{#1}}}

\newcommand\eight{\fontsize{8pt}{8pt}\selectfont}

\newcommand{\circH}{\Theta}

\begin{document}

\title{The Cost-Accuracy Trade-Off In Operator Learning\\ With Neural Networks}


\author[1]{
Maarten V. de Hoop\thanks{ 
  {\tt mdehoop@rice.edu}.
  }
}
\author[2]{
Daniel Zhengyu Huang\thanks{ 
  {\tt dzhuang@caltech.edu}.
  }
}
\author[2]{
Elizabeth Qian\thanks{ 
  {\tt eqian@caltech.edu}.
  }
}
\author[2]{
Andrew M. Stuart\thanks{ 
  {\tt astuart@caltech.edu}.
  }
}

\affil[1]{Rice University, Houston, TX, USA.}
\affil[2]{California Institute of Technology, Pasadena, CA, USA.}

\begin{abstract}
The term `surrogate modeling' in computational science and engineering refers to the development of computationally efficient approximations for expensive simulations, such as those arising from numerical solution of partial differential equations (PDEs). Surrogate modeling is an enabling methodology for many-query computations in science and engineering, which include iterative methods in optimization and sampling methods in uncertainty quantification. Over the last few years, several approaches to surrogate modeling for PDEs using neural networks have emerged, motivated by successes in using neural networks to approximate nonlinear maps in other areas.
In principle, the relative merits of these different
approaches can be evaluated by understanding, for each one, the
cost required to achieve a given level of accuracy. However, the absence of a complete theory of approximation error for these approaches makes it difficult to assess this cost-accuracy trade-off. The purpose of the paper is to provide a careful
numerical study of this issue, comparing a variety of different
neural network architectures for operator approximation across
a range of problems arising from PDE models in continuum mechanics.
\end{abstract}

\keywordone{Computational partial differential equations,}
\keywordtwo{Surrogate modeling,}
\keywordthree{Operator approximation,}
\keywordfour{Neural networks,}
\keywordfive{Computational complexity,}

\maketitle

\section{Introduction}
\label{sec:I}
In many problems in computational partial differential equations (PDEs) the fundamental driver in deciding which approximation methodology to employ is the shape of the cost-accuracy curve: this determines what computational resources are required to achieve a desired level of accuracy, a measure of
computational complexity. On this basis some methods may be shown to clearly 
outperform others, guiding computational practice. In the numerical analysis of PDEs there is a deep literature addressing this issue.
This literature comprises two main components: (i) an analysis of the error as a function of the resolution of the finite dimensional approximation~\cite{jovanovic2013analysis,strikwerda2004finite,johnson2012numerical,canuto2012spectral,trefethen2000spectral}; and (ii) analysis of the cost of running the model, at a given level of finite-dimensional resolution, often dominated by 
matrix inversion and/or matrix-vector multiplies  \cite{demmel1997applied,trefethen1997numerical}
and/or by time-stepping and iteration count for nonlinear solvers. Theoretical results in (i) and (ii) may be combined to determine the cost-accuracy curve for different methods and thereby inform the choice of method for a given problem.
For certain classes of equations, multi-resolution methods have emerged which are near optimal in terms of minimizing cost for a given error~\cite{hackbusch2013multi,owhadi2019operator}.

Data-driven approximation of mappings/operators between function spaces provides a way to learn cheap-to-evaluate surrogates
which can bypass the need for \eq{employing PDE solvers}, after an initial
training phase in which data are generated. These surrogates then enable efficient many-query analyses of PDE-based problems in computational science and engineering.
However, the theory for data-driven approximations 
is in its infancy and cost-accuracy curves are not analytically understood. The goal of this work is to provide a numerical study of the cost accuracy trade-off, for a range of operator neural network architectures, including PCA-based neural networks~(PCA-Net)~\cite{hesthaven2018non,bhattacharya2021model}, DeepONet~\cite{lu2021learning,lu2021comprehensive}, pointwise evaluation~(PARA-Net, defined in this paper), and the Fourier neural operator~(FNO)~\cite{li2020fourier,kovachki2021neural}. 
The numerical studies are conducted on four test problems: (1) the two-dimensional incompressible Navier-Stokes equation, (2) the Helmholtz equation, (3) a structural mechanics test problem, and (4) the linear advection equation.

There are four sources of error in these operator learning problems: a) discretization of the input and output spaces; b) parameterization of the operator approximators; c) finite data volume; d) the optimizer used in training. 
In this paper we concentrate on b) and c) and study the cost-accuracy trade-off in relation to data volume and number of parameters in the neural network. The reason for not studying a) in this work is that,
if properly designed, operator approximators have the property of discretization invariance, meaning that they are defined to act between function spaces and training of parameters for one discretization can therefore  be used for other discretizations \cite{bhattacharya2021model,li2020fourier,kovachki2021neural}; in this setting a single set of parameters will provide good approximations for all resolutions for which the discretization error is small enough.
\eq{As for the role of the optimizer d), while there exists numerical evidence that stochastic gradient descent methods can be effective in driving the loss function (close) to its global minimum~\cite{goodfellow2016deep,hardt2016train,zhang2021understanding}, this work is far from being theoretically well-understood and, furthermore, is not in the context of operator learning and partial differential equations. There are also other optimization approaches, for example using second-order information \cite{nocedal1999numerical,bottou2018optimization}, or ensemble methods \cite{haber2018never,kovachki2019ensemble}, that may produce different results. However, in order to limit the scope of our numerical studies, we employ stochastic gradient descent using fixed standard choices of the optimization hyperparameters for all test cases. When there is evidence that the optimization itself limits the accuracy achieved, we will highlight this in our discussion.
With this caveat, we primarily focus our study on b) and c) and extract clear signals from numerical experiments, laying foundation for future studies which delve into the interactions with a) and d).}
Our numerical experiments will disentangle the roles of errors caused by b) and by c).
The seminal work of Giles \cite{giles2008multilevel,giles2015multilevel} on multilevel Monte Carlo methods demonstrates that a theoretical understanding of errors incurred through the interaction of finite sampling and finite dimensional approximation leads to highly efficient methods which use different sample sizes for different finite
dimensional approximations of expectations. Future analysis studying the interaction between the sources of error arising from b) and c) would be very valuable in the field of operator approximation and could lead to improved complexity results, 
\eq{in the spirit of this work of Giles. Furthermore, although we
downplay the effect of discretization error a) in this paper,
future analysis studying the interaction between the sources of error arising from a), b) and c) could be also} very valuable in minimizing the cost of the objective evaluation during training.

\subsection{Literature Review}
\paragraph{Many-Query Motivation}
Many computational tasks arising in science and engineering require repeatedly evaluating the output of an expensive computational model, which generally involve solving parametric partial differential equations.
Examples include design optimization~\cite{economon2016su2,martins2013multidisciplinary,bendsoe2003topology,boncoraglio2021active}, Bayesian calibration and inference~\cite{stuart2010inverse,martin2012stochastic,huang2022iterated,huang2022efficient,morzfeld2019localization,cao2021bayesian,herrmann2020deep}, and multiscale computation~\cite{weinan2011principles,fish1997computational,feyel2000fe2,liu2022learning,kovachki2021multiscale,chen2021exponential,chen2021exponentially}.
In some settings, this computational model may be viewed as a black box mapping functions to functions, making the development of efficient data-driven surrogate models highly desirable.
Furthermore, there are increasingly complicated application domains which are data-rich, but for which  complexity limits first principles modeling;
thus model discovery from data is necessitated.
Cyber-physical systems present a wide-class of such problems \cite{song2016cyber}. Data-driven modeling has great potential to impact such domains.

\paragraph{Data-driven Surrogate Models}
A variety of specific techniques for 
data-driven surrogate modeling
have been developed and described in the literature, including the Koopman operator-based models~\cite{schmid2010dynamic,rowley2009spectral,tu2013dynamic,taira2017modal,mezic2005spectral,mezic2013analysis}, Gaussian process~(GP) based model emulators~\cite{kennedy2001bayesian,cleary2021calibrate,chen2021solving}, and data-driven projection-based reduced models~\cite{devore2014theoretical, carlberg2013gnat,peherstorfer2016data,binev2017data,cohen2020state,qian2020lift,qian2021reduced,grimberg2021mesh,gosea2021data}.
Some of these methods are non-intrusive, or can be extended to non-intrusive instantiations, which can be constructed without prior knowledge of an underlying mechanistic model. The motivation for such methods is to find fast but accurate
replacements for expensive computational tasks which need to be repeatedly executed.

\paragraph{Neural Network Based Surrogate Models}
Most of the aforementioned surrogate models have at their core a linear model represented as combinations of appropriate basis functions. Composing such models with pointwise nonlinearities in layers leads to neural network based surrogate models.  The introduced nonlinearity in the surrogate models does not lead to significant increase in evaluation cost, since it is
pointwise, but can significantly improve expressivity and prediction accuracy. As a consequence, neural network based surrogate models are being explored in different fields in science and engineering. They are used for data-driven modeling and scientific discovery, including turbulence modeling~\cite{ling2016reynolds,wang2017physics,wu2018physics,duraisamy2019turbulence}, material modeling~\cite{ghaboussi1998autoprogressive,zhang2018deep,avery2021computationally,liu2022learning,liu2021a}, quantum mechanical modeling of materials \cite{zhang2018deep,lu202186} and the design and prediction of protein structures~\cite{jumper2021highly}. 
In addition, neural networks have been used to augment conventional data-driven surrogate modeling frameworks. For example, proper orthogonal decomposition is coupled with neural networks to speed up the online stage of reduced models in~\cite{hesthaven2018non,swischuk2018physics}, nonlinear manifolds are introduced and learned by convolutional autoencoders to extend projection-based model reduction beyond approximation in linear subspaces in~\cite{lee2020model}, the approximation of the Koopman operator is enhanced by neural networks for the control of unsteady fluid flows
~\cite{li2017extended,morton2018deep}, relationships between GPs and deep neural networks \eq{are} explored in~\cite{lee2017deep}, and Bayesian neural networks for uncertainty quantification are investigated in ~\cite{zhu2018bayesian,yang2021b}.

\paragraph{Neural Networks for Solving Partial Differential Equations}
Solving partial differential equations with neural network-based approximations provides an alternative to traditional approximation approaches (such as finite difference, finite element or spectral methods) by introducing new forms of finite dimensional approximation spaces.
Potential advantages include the fact that neural network-based approximations do not require meshes, and also offer the
possibility of exploiting nonlinear approximation spaces, both of
which may be of great benefit in the numerical solution of high-dimensional PDEs \cite{han2017deep,sirignano2018dgm,han2018solving}, and in solving PDEs on complex geometries~\cite{berg2018unified,raissi2019physics,goswami2020transfer}. 
Moreover, well-developed neural network libraries (i.e., TensorFlow~\cite{tensorflow2015-whitepaper} and PyTorch~\cite{paszke2019pytorch}), coupled
with GPU architecture, make the solution process both easy-to-use and efficient. The basic methodology in PINNs~\cite{raissi2019physics} introduces a loss function defined by the equation, and may be extended to inverse problems through the additional incorporation of data-related penalty terms~\cite{raissi2020hidden,sun2020physics,mao2020physics}. The related gradient can generally be computed automatically by these well-developed neural network libraries, and hence the optimization or inverse solving procedure are also both easy-to-use and efficient.
Theoretical underpinnings of this framework are presented in~\cite{he2018relu,weinan2019barron,shin2020convergence,daubechies2021nonlinear} and approaches to improving the performance with locally adaptive activation functions are studied in ~\cite{jagtap2020adaptive,jagtap2020locally,hu2021extended,jagtap2022deep}.

\paragraph{Neural Networks for Operator Learning}
In many of the applications cited in the preceding paragraphs, in both surrogate modeling
and in model discovery, the core task is the mapping of one function into another. Much of
the preceding work on neural networks for PDEs can be naively extended to neural network methods which would need to be re-trained for each different input instance; some
of the reduced order modeling \eq{work} proceeds by learning mappings on finite dimensional spaces without explicitly conceptualizing or recognizing the function space mapping at its core.
The focus of the present work is on infinite-dimensional formulations for learning mappings between function spaces, and we refer to this as operator learning. There has been considerable
activity in this area in the last five years\eq{:} starting with the paper \cite{zhu2018bayesian}, focused
on surrogate modeling for uncertainty quantification in subsurface flow using deep 
autoencoder networks\eq{;} \cite{khoo2021solving} who use a convolutional neural  network to minimize parameterized variational problems\eq{;} DeepONet~\cite{lu2021learning} with architecture
comprising two sub-networks used to represent the operator by enforcing separation of input 
and output functions into the branch and trunk networks\eq{;} PCA-Net~\cite{hesthaven2018non,bhattacharya2021model} which uses neural networks to map between PCA coefficients representing input and output functions, neural networks based on kernel integral operators \cite{kovachki2021neural}, and the FNO~\cite{li2020fourier,wen2021u,li2022fourier} which parameterizes the integral kernel directly in Fourier space. \eq{Finally, we mention
\cite{gin2020deepgreen,boulle2022learning} which study the learning of Green’s functions for nonlinear boundary value problems,} \cite{khoo2021solving,fan2020solving,fan2019solving,khoo2019switchnet} which propose  novel neural network architectures for solving wave and elliptic equations, and \cite{kissas2022learning} which employs the attention mechanism for operator learning.
A key aspect of the work in \cite{bhattacharya2021model,kovachki2021neural} is that the methodology has been designed to be robust to the resolution of the finite-dimensionalizations
of the input and output functions; this means that the computational investment in learning
a model at one resolution, or for one particular discretization, may be transferred
to other resolutions or discretizations. \eq{There is an existing body of published work
which compares resolution invariant methods such as PCA-NET, DeepONet and FNO,
with methods where the choice and training of parameters is linked to the
grid resolution, such as U-NET and other grid-adapted approaches \cite{zhu2018bayesian};
see \cite{li2020fourier,bhattacharya2021model} for such comparisons. These comparative studies
demonstrate the value of using resolution invariant methods.}

Universal approximation theorems, 
stating that that neural networks can accurately approximate  wide classes of nonlinear 
continuous operators are starting to emerge
to underpin these methodologies~\cite{chen1995universal,lu2021learning,kovachki2021neural}.
The paper \cite{lu2021comprehensive} undertakes a careful comparison of the relative merits of
DeepONet-based and FNO-based methodologies for operator learning. The paper introduces
variants on the basic methods and the results show the impressive flexibility of DeepONet and its
ability to solve problems from a variety of physical applications and a variety of complex geometric
domains, and suggest directions in which the FNO will benefit from further innovations. The
specific numerical results presented, in which a fixed number of parameters is chosen for
each method, and the error compared, show smaller errors for DeepONet than for FNO in many
examples; however similarly designed experiments on different problems and with different
choices of the fixed parameters reach the opposite conclusion about the ordering of
the errors incurred by DeepONet and FNO \cite{kovachki2021neural,kissas2022learning}. Thus, in terms
of evaluating the relative merits of different surrogates,
the papers \cite{lu2021comprehensive,kovachki2021neural,kissas2022learning} come up short,
because they do not  study
dependence of evaluation cost on the approximation error achieved. This question
is the primary focus of our paper.

\subsection{Our Contributions}
Neural network surrogates for operators require a significant volume of training data to ensure reasonable predictive power. The number of parameters needs to increase with data volume, and when the parameterization of the neural network architecture becomes more complicated, the evaluation cost increases. Quantifying the resulting cost-accuracy trade-off involves interaction between statistical error from data sampling and approximation error from parameterization; in particular the issue of how to choose the number of parameters, given
the amount of available data, is fundamental to success and efficiency of the methodology. Understanding the trade-off between cost and accuracy, for different neural operators, and in different parameter/data regimes, is
thus of central importance in guiding how this field develops; it forms the focus of this work. Because theory in this arena is currently limited, our study is purely computational.
Our belief is that careful computational studies will provide impetus for the development of theoretical understanding of the issues they explore.  Our main contributions are as follows:

\begin{itemize} 

\item We give a unified presentation of a variety of different operator approximators: three have appeared in the recent
literature,  PCA-Net, DeepONet and FNO, whilst a fourth, a lifting of finite dimensional
neural networks via pointwise evaluation,
which we label PARA-Net, is formalized in this paper.

\item We numerically study these operator approximators in the
context of four model problems: 
    i) mapping forcing to solution in a 2D incompressible Newtonian fluid; ii) mapping wavespeed to disturbance field in the Helmholtz equation; iii) mapping applied tension load to stress field in an elastic solid; and iv) mapping initial condition to solution at time one in the advection equation.
     
    \item We study error as a function of data volume (for fixed
    number of parameters), error as a function of number of 
    parameters (for fixed data volume), and evaluation cost as
    a function of error,  for all four test problems. \eq{These experiments 
    quantify anticipated differences in the cost-accuracy trade-off between problems with smooth outputs (i,ii) and those with discontinuous outputs (iii,iv). The experiments also demonstrate data limiting effects and the potential for over-fitting in some methods.}
    
    \item We show that the PARA-Net approach, whilst being a natural generalization of
    neural networks between Euclidean spaces, is not competitive with the other three methods considered, for the problems we consider, as measured by evaluation cost per unit accuracy. PCA-Net, DeepONet and FNO all exhibit desirable behavior, and are all clearly viable methodologies for certain problems in certain regimes; in contrast PARA-Net 
    is both consistently more expensive whilst also being particularly prone to over-fitting, except for
    problem (iv); as a result it does not appear to be a useful general approach. We include PARA-Net because doing so highlights the drawbacks of not conceptualizing operator approximation as learning a function to function mapping.
    
    \item For a number of test problems and operator
    approximations, we provide explicit instances of the 
    test cases which lead to the median test error
    and to the largest test error, yielding insight 
    into failure modes of the learning procedures,
    especially for the non-smooth problems (iii,iv).
    
    \item We numerically study the range space 
    of the DeepONet approximator, contrasting it with PCA-Net, demonstrating the pros and cons of its basic implementation, and motivating the PCA modification of the basic DeepONet method that is introduced in \cite{lu2021comprehensive}.
\end{itemize}
\vspace{0.12in}

Our experiments do not provide a definitive answer as to which operator approximator is best for any given class of problems. Rather we focus our attention on the following {\em complexity} question: how does 
evaluation cost scale with accuracy for these methods? Our experiments demonstrate that the answer to this question depends on the problem. In a similar vein, the results in \cite{lu2021comprehensive,kovachki2021neural,kissas2022learning} also indicate that preferred method with respect to error at a fixed parametrization level is problem dependent. However, since \cite{lu2021comprehensive,kovachki2021neural,kissas2022learning} consider only fixed parameterizations, and hence fixed cost, they do not shed light on complexity. 
Our experiments demonstrate the merits of considering the complexity question, and suggest the need for theory to underpin empirical findings. We believe conclusions as to preferred methodology for any given problem, even with respect to the complexity measure studied here, are premature; our numerical
results simply focus on the need for theory which addresses the complexity question.

The code and dataset are accessible online at  \url{https://github.com/Zhengyu-Huang/Operator-Learning}.

\section{Problem Formulation}
\label{sec:PF}
Consider the following Banach spaces of functions $\cU, \cV$:
\begin{align*}
    \cU = \{u: D_u \rightarrow \R^{d_i}\}, \textrm{ where }  D_u \subseteq \R^{d_x}\eq{,} \\
    \cV = \{v: D_v \rightarrow \R^{d_o}\}, \textrm{ where }  D_v \subseteq \R^{d_y}\eq{.}
\end{align*}
Our goal is to determine operators $\Psid : \cU \rightarrow \cV$ from (samples from) the probability measure $\bigl({\rm Id},\Psid\bigr)^{\#} \mu$, where $\mu$ is supported on $\cU$ and equipped with the Borel $\sigma$-algebra, and the push-forward $\left({\rm Id},\Psid\right)^{\#} \mu$ is supported on $\cU \times \cV$, also
equipped with the Borel $\sigma$-algebra. For simplicity
we will assume that  $D_u$ and $D_v$ are closed and bounded.

We approximate $\Psid$ by different classes of parametric operators, each of which characterizes a different neural network architecture. To this end, let $\circH \subseteq \R^p$ denote the space of parameters of the neural network and consider a family of functions from $\cU$ into $\cV$
defined by
$\Psi: \cU \times \Theta \mapsto \cV.$
In the idealized case
of infinite data, the parameters $\theta\in\Theta$ are chosen to be $\theta^\star$, the minimizer of the risk defined by: 
\begin{align*}
\text{Risk:}\,\, &                        \cR_{\infty}(\theta):=\E^{u\sim\mu}\|\Psi^\dagger(u)-\Psi(u;\theta)\|_{\cV}^2=\|\Psi^\dagger-\Psi\|^2_{L^2_{\mu}(\cU,\cV)},
\end{align*}
where the last expression is defined with integration over
$\cU$ in the Bochner sense. In practice we only have access to
samples $\{u_n\}_{n=1}^N$ from $\mu$, assumed to be 
i.i.d.\ and defining the resulting empirical measure $\mu^N$, together with $\{\Psi^\dagger(u_n)\}_{n=1}^N$.
\eq{This enables us to define the empirical risk:}
\begin{align*}
\text{Empirical Risk:}\,\, &                        \cR_{N}(\infty)(\theta):=\E^{u\sim\mu^N}\|\Psi^\dagger(u)-\Psi(u;\theta)\|_{\cV}^2=\frac{1}{N}\sum_{n=1}^N\|\Psi^\dagger(u_n)-\Psi(u_n;\theta)\|^2_{\cV}.
\end{align*}
Parameter $\theta^{\star,N}$ is the
value of $\theta$ which minimizes $\cR_{N}(\infty)(\theta).$
In practice, in view of the non-convex nature of the optimization
over $\theta$, we may only have access to an approximation of
$\theta^{\star,N}$.

\begin{remark}
Within both the risk and the empirical risk various modifications may be relevant; for example
$\lVert \cdot \rVert_{\cV}$ may be replaced by the norm in any space into which $\cV$ is continuously embedded.
When studying error as a function of data volume, parameter dimension, and evaluation cost,
we will also use the following relative error measure (approximated empirically):

\begin{align}
    \E^{\mu} \left(\frac{\lVert\Psi^\dagger(u) - \Psi(u; \theta^{\star,N}) \rVert_{\cV}}{\lVert\Psi^\dagger(u) \rVert_\cV}\right).\label{eq: relative error measure}
\end{align}

\end{remark}

\begin{remark}
Many of our examples concern the case where $D_u=D_v=D$, $d_x=d_y=d$.
Furthermore, many of our examples also concern the case where the vector-valued functions in $\cU$
and $\cV$ are  in fact scalar-valued so that
the input and output dimensions are $1$ : $d_i = d_o = 1$.
However, since the methodology applies outside these simpler settings, we expose basic ideas at this greater level of generality.
\end{remark}

\section{Neural Networks}
\label{sec:NN}

In this section we describe the various neural networks that
we will compare in this study. We establish preliminary notation
in \Cref{ssec:P}. In  \Cref{ssec:PCA,ssec:DEEPO,ssec:PARA,ssec:NEU}, we introduce four classes of neural network, highlighting similarities
and differences. 

\subsection{Preliminaries}
\label{ssec:P}

Let $\cH$ denote a Hilbert space with inner product $\left<\cdot,\cdot\right>$ and induced norm $\|\cdot\|$.
We then consider the Gelfand triple $\cU \subseteq \cH \subseteq \cU^{*}$
where the embedding of $\cU$ in $\cH$ is continuous and $\cU^*$ denotes
the dual space of $\cU$ -- the space of linear functionals on $\cU.$
Given $L_k \in \cU^{*},\,k=1,\cdots d_u$, the function $u$ in $\cU$ is partially characterized by the finite dimensional vector $Lu:=\{L_k u\}_{k=1}^{d_u}$. This will be used to characterize inputs to $\Psi$ in the formulations of the three neural network architectures in \Cref{ssec:PCA,ssec:DEEPO,ssec:PARA}.

When $\cU,\cV$ are themselves Hilbert spaces we may compute the mean and covariance operator with respect to $\mu$ and $(\Psid)^{\#} \mu$, respectively, assuming first and second moments exist. The eigen-decomposition of the covariance  operators (PCA) gives rise to orthonormal bases 
\eq{$\{\phi_j\}_{j\in\N}$ in $\cU$ 
 and $\{\psi_j\}_{j\in\N}$ in $\cV.$
Note that
\begin{equation*}
\phi_j : D_u \rightarrow \R^{d_i},\quad  j\in\N;\quad\quad
\psi_j : D_v \rightarrow \R^{d_o}, \quad j\in \N.
\end{equation*}
}

\subsection{PCA-Net}
\label{ssec:PCA}
Here it is simplest to think of $\cH=\cU=L^2(D_u;\R^{d_i}).$
Truncate the PCA basis of $\mu$ to the first $d_u$ modes, $\{\phi_j\}_{j=1}^{j=d_u}$, and truncate the PCA basis of $\Psi^{\#} \mu$ to the first $d_v$ modes $\{\psi_j\}_{j=1}^{j=d_v}$.
Then define $L_k \in \cU^{*}$ by $L_k u = \left<\phi_k, u\right>,\, k = 1,\,\cdots d_u$. Introduce function
$\alpha: \R^{d_u} \times \circH \rightarrow \R^{d_v}$ defined componentwise by
$$\alpha_j : \R^{d_u} \times \circH \rightarrow \R \qquad j = 1,2,\cdots,\ d_v$$ 
which maps from the PCA coefficients $\{L_k u\}_{k=1}^{d_u}$ of $\mu$ to the PCA coefficients $\{\alpha_j\}_{j=1}^{d_v}$ of $\Psi^{\#} \mu$, and is parameterized by $\theta \in \Theta.$ 
Here function $\alpha$ denotes a finite dimensional neural network.

The approximating operator between $\cU$ and $\cV$ is then defined by
\begin{align*}
    \Psi_{PCA}(u;\theta)(y) = \sum_{j=1}^{d_v}\alpha_j(L u;\theta) \psi_j(y) =  \alpha(L u;\theta )^T \psi(y) \qquad \forall u\in U \qquad y\in D_{v}
\end{align*}
where $L u = \{L_k u\}_{k=1}^{d_u}$ and $\psi(y) = \{\psi_j(y)\}_{k=1}^{d_v}$.
We note that each $\alpha_j$ is $\R$-valued, so $\alpha$ may be
viewed as a column vector in $\R^{d_v}$, whilst each $\psi_j$ is
$\R^{d_o}$-valued, so $\psi$ may be viewed as a matrix in $\R^{d_v} \times \R^{d_o}.$ Thus $\Psi_{PCA}(u; \theta, \mu)(y)$ is $\R^{d_o}$-valued.

\begin{remark}
Note the following facts concerning this architecture which place it
slightly outside the standard training framework, using empirical risk, outlined above.
\begin{itemize}
\item \eq{PCA is defined on a Hilbert space, given measure $\mu$; in practice it is
implemented on finite-dimensional approximations of a Hilbert space, using samples
from measure $\mu.$}
    \item  The linear functionals $\{L_k\}$ in $\cU^*$ and the functions $\{\psi_j\}$ in $\cV$, needed to define the architecture of 
$\Psi_{PCA}$, are pre-computed using PCA and decoupled
from the training of the neural network. Thus preliminary calculations,
\eq{using PCA on the dataset}, are conducted to define the risk or empirical risk.
\end{itemize}
\end{remark}





\subsection{DeepONet}
\label{ssec:DEEPO}
\eq{Several works in recent years use the terms `DeepONet' or `deep operator network' to describe related, but different neural network architectures used to approximate infinite-dimensional operators \cite{lu2019deeponet,lu2021learning,goswami2020transfer,goswami2022physics,sharma2021application}. These architectures have in common the use of two separate networks, called the \emph{trunk} and \emph{branch} networks, to map input functions to output functions,
generalizing the `shallowONet' architecture first proposed in \cite{chen1995universal}.
The trunk network is used learn the output space representation; 
the branch network learns the input to output solution mapping in the output space 
spanned by the trunk network.}

\eq{In this paper we employ a specific variant of the DeepONet branch/trunk
architecture, chosen in order to faciliate comparison with PCA-Net; 
specifically we consider a version of DeepONet in which we use the same linear 
functionals in $\cU$ as in PCA-Net as input to the branch network. However, whilst PCA-Net
uses fixed PCA basis functions in the output space, DeepONet \emph{learns} the
representation in the output space as the trunk network. We note that
the variant on DeepONet appearing in~\cite[Problem 5]{lu2021learning} also uses
a Karhunen-Loeve expansion, which is equivalent to PCA,
as input, but uses the information
in a different way; in particular the coefficients of the Karhunen-Loeve expansion, which in our notation are analogous to $Lu$,
are used in~\cite{lu2021learning} as input to the trunk network, rather than the branch network.}

\eq{In summary the key difference between our implementations of DeepONet and PCA-Net is that the functions $\{\psi_j\}$ in the output space are given in DeepONet by a neural network, rather than using PCA on the output space $\cV$ as in PCA-Net; for DeepONet they are learned in the training phase, rather than being precomputed from the data as in  PCA-Net.  However, our implementations of DeepONet and PCA-Net
both take linear functionals on the input PCA basis in $\cU$
as input to the neural networks labelled $\alpha$. 
We emphasize that the original formulation of the DeepONet used pointwise evaluations of the input function in $\cV$ as inputs to the neural network, rather than PCA coefficients.
A brief discussion of pointwise evaluations, written in our general framework, 
may be found in Appendix \ref{sec:pw}. It is possible, indeed likely, that for some
problems, our implementation of both DeepONet and PCA-Net could show improved performance
by working with pointwise evaluation functionals as input; however,
our numerical results are not focused on this question, 
but rather on how the two methods perform in relation to their different representations
of the output space.}

The operator between $\cU$ and $\cV$ is defined as follows:
\begin{align*}
    \Psi_{DEEP}(u;\theta)(y) = \sum_{j=1}^{d_v}\alpha_j(L u;\theta_\alpha) \psi_j(y;\theta_{\psi}) = \alpha(L u;\theta_{\alpha} )^T \psi(y;\theta_{\psi}), \qquad \forall u\in U, \qquad y\in D_{v},
\end{align*}
where the functions $\alpha_j : \R^{d_u} \times \circH \rightarrow \R$ have
the same structure as in PCA-Net, and are collectively referred to as the ``branch'' of the DeepONet neural network. The functions $\{\psi_j\}$ are collectively referred to as the ``trunk'' of the DeepONet architecture, and are defined by neural networks of the form
\begin{align*}
    \psi_j : D_v \times \circH \rightarrow \R^{d_o}, \quad j = 1,\cdots , d_v.
\end{align*}
We denote by $\theta$ the collection of the hyperparameters $\theta_\alpha$ and $\theta_\psi$ appearing in $\alpha$ and $\psi$. These are computed from minimizing the (empirical) risk over $\theta \in \circH.$ No parameters in $\circH$ are shared  between $\alpha$ and $\psi$, but \eq{their optimal choice} is coupled 
through the empirical risk minimization.

\begin{remark} Our choice of notation highlights
similarities between PCA-Net and DeepONet. However
features that distinguish DeepONet from PCA-Net include:
\begin{itemize}
    \item For DeepONet, both the trunk, $\{\psi\}$, and the branch, $\{\alpha\}$ are neural networks, whereas for PCA-Net only the
    branch is. In particular for DeepONet $\psi(\cdot;\theta)$ is  
    selected via optimization during training, 
    in contrast to PCA-Net for which $\psi$ is defined by the data, using PCA,
    prior to parameter selection, via optimization, for the neural network.
    
    \item In DeepONet,  $L$ may not depend 
    on $\mu$ if pointwise evaluations are used (see Appendix \ref{sec:pw}.)

\item For both PCA-Net and DeepONet, the norm $\lVert \cdot \rVert_{\cV}$ appearing in the empirical risk, 
or weaker norm in a space into which V is continuously embedded, is approximated using pointwise evaluations at a set of points $\{y_l\}$. Indeed, in the original work underpinning DeepONet \cite{chen1995universal} a set of fixed pointwise evaluations are part of
the definition of the architecture and appear in the definition of the empirical risk used in training. We prefer to employ an operator perspective on the method, and in particular to define the empirical risk through a norm on $\cV$, enabling comparison with other
neural networks formulated as operator approximators.
\end{itemize}
\end{remark}

\subsection{PARA-Net}
\label{ssec:PARA}

As we do for PCA-Net, we truncate the PCA basis of $\mu$, $\{\phi_j\}_{j=1}^{d_u}$, using the first $d_u$ modes, and define 
$L_k \in \cU^{*}$ by $L_k u = \left<\phi_k, u\right>,\, k = 1,\,\cdots d_u$. 
Again it is simplest to think of $\cH=\cU=L^2(D_u;\R^{d_i}).$
We introduce the real-valued neural network
\begin{align*}
    \psi: \R^{d_u} \times D_v \times \Theta \rightarrow \R^{d_o}
\end{align*}
and then define $\Psi_{PARA}: \cU \times \circH \rightarrow \cV$  by 
\begin{align*}
    \Psi_{PARA}(u;\theta)(y) = \psi(Lu, y;\theta).
\end{align*}


\begin{remark}
Note that $\psi$ is a standard neural network
between Euclidean spaces. By defining a collection of
linear functionals $L$ in $\cU^*$, evaluating this finite
dimensional neural network at input $(Lu,y)$, and varying
over $u$ in $\cU$ and $y$ in $D_v$ we create an operator.
In this sense it is a natural method. However, as we will show, it is not competitive with the other methods presented
here in terms of the cost-accuracy trade-off for the problems we consider.
We include it because it is a natural generalization  but its poor performance
serves as a motivation to think beyond the confines of
standard finite dimensional Euclidean neural networks and to exploit structure such as PCA bases, branch-trunk decomposition or (next subsection) Fourier representation.

Recall that the graph of $\Psi$ is the set 
$\{\Psi(u), u \in \cU\}.$
Both PCA-Net and DeepONet give rise to linear approximation spaces for the graph of $\Psi$ in that, for all $u \in \cU$, the approximation
of $\Psi(u) \in \cV$ lies in the same finite dimensional subspace of $\cV$, independently of $u$, 
defined by the span of the $\{\psi_j\}.$ 
In contrast PARA-Net gives rise to a 
an approximation space in $\cV$ which depends nonlinearly on $u.$
This is a potential advantage but, as we will show for the
problems we consider here, this advantage is outweighed by the computational cost of the pointwise evaluations required to construct full output function reconstructions. \eq{This challenge is exacerbated by multiple spatial dimensions, since the number of grid points generally scales exponentially with spatial dimension.} However it is conceivable
that, for some operator approximation problems not considered here,
the trade-off between nonlinear approximation and the cost 
resulting from repeated pointwise evaluation results in PARA-Net
being competitive. The potential benefits of nonlinear approximation spaces
are discussed in \cite{devore1998nonlinear}.
\end{remark}

\subsection{Fourier Neural Operator (FNO)}
\label{ssec:NEU}

For simplicity, we consider the setting where $D_u = D_v = D= [0,1]^d$ so that $d_x=d_y = d$.
In this setting it is simplest to think of $\cH=L^2(D;\R^{d_i})$
(or a generalization to impose, for example, the divergence-free
condition)
and $\cU=\cV=C_{\rm per}(D;\R^{d_i})$, the set of continuous periodic functions on the unit cube; alternatively $\cU$ may be an RKHS,
such as a Sobolev space of periodic functions of fractional order greater than $d/2$,
which is continuously embedded into $C_{\rm per}(D;\R^{d_i}).$

Let $R, Q$ denote standard finite dimensional neural networks
\begin{align*}
    R: \R^{d_i} \times \circH \rightarrow \R^{d_f}\\
    Q: \R^{d_f} \times \circH \rightarrow \R^{d_o},
\end{align*}
where $d_f$ is the number of channels\footnote{Sometimes
also referred to as features; however we avoid this terminology  
as the terminology random features are used to describe a different concept \cite{nelsen2021random} and we wish to avoid confusion.}
used in FNO-NET, typically larger than $d_i$ or $d_o.$
We use $R$ and $Q$ to define operators $\cR$ and $\cQ$ by pointwise evaluation:
\begin{equation}
    \label{eq: FNO lift downsample defs}
    \begin{aligned}
        (\cR u)(x, \theta_R) &= R(u(x), \theta_R)\\
        (\cQ u)(x, \theta_Q) &= Q(u(x), \theta_Q). 
    \end{aligned}
\end{equation}
We say the $\cR$ operator \textit{lifts} the input to the
channels and that the $\cQ$ operator \textit{projects} 
the channels to the output.
Note that $\cR u$ is well-defined
in $C_{\rm per}(D;\R^{d_f})$, the space of periodic continuous functions of dimension $d_f$, if $u \in \cU$ and $\cU$ is continuously embedded into $C_{\rm per}(D;\R^{d_i}).$

We now introduce notation for the $l$-th \textit{Fourier Neural Layer} (FNL): 
\begin{align}
    \cL_{l}(v)(x, \theta) = \sigma\big(W_l v(x) + (\cK v)(x;\gamma_l)\bigr),
    \label{eq: FNL definition}
\end{align}
where  $\sigma$ is an activation function, applied pointwise w.r.t. $x \in D$; $W_l \in \R^{d_f \times d_f}$ is a matrix applied pointwise to $v(x)$;  $\cK$ is a parameterized non-local operator, 
for which a variety of forms are commonly used \cite{kovachki2021neural}. In this paper
we exclusively use the Fourier Neural Operator form (FNO):
\begin{align*}
    \text{FNO}: \quad (\cK v)(x; \gamma) &= \cF^{-1} (P(\gamma) (\cF v))(x).
\end{align*}
Here 
$\cF$ denotes the Fourier transform of a periodic function $v: D \rightarrow \R^{d_f}$, so that $\cF v:  \Z^{d} \rightarrow \C^{d_f}$, and $\cF^{-1}$ denotes its inverse,
and for each point in the Fourier domain $\Z^{d}$, the action of $P(\gamma)$ is as an
element of $\C^{d_f \times d_f}$.
In this paper the FNO is implemented using a Fast Fourier Transform on a uniform lattice of pointwise evaluations of the functions $v$ and $P(\gamma)\cF(v)$. When evaluating $\cF$, only the first $k_{max}$ modes are kept.
In this setting $P$ reduces to a complex-valued ($k_{max} \times d_f \times d_f$)-tensor, 
which is applied as an element of $\C^{d_f \times d_f}$  on each of the $k_{max}$ Fourier modes. Then
\begin{align*}
    \Psi_{FNO}(u;\theta) = \cQ\circ \cL_{L}\circ \cdots \cL_{2}\circ \cL_{1}\circ \cR(u).
    \end{align*}
Here we denote by $\theta$ the collection of the hyperparameters $\theta_R$ and $\theta_Q$ and $\{W_l, \gamma_l\}$ for each Fourier neural layer $l$. As in the case of DeepONet, none of these hyperparameters are shared, but their choice is coupled through the empirical risk minimization.
Note that, for simplicity, the layer width $d_f$ is fixed, but this is not necessary.
Furthermore, all numerical experiments reported here are conducted with $Q$ and $R$ being affine. In contrast, the experiments in \cite{li2020fourier,kovachki2021universal} 
use linear $R$ but nonlinear $Q.$

\eq{Finally, we note that in the current implementation of FNO, the number of retained Fourier modes, $k_{max}$, scales exponentially with the spatial dimension of the problem. In contrast, our implementations of PCA-Net and DeepONet do not have parameters that explicitly depend on the spatial dimension of the problem. The FNO scaling with spatial dimension may pose a computational hurdle in three spatial dimensions; however we note that this may be overcome, e.g.\ by employing wavelets in higher dimensions. }

\section{Numerical Studies}
\label{sec:NUM}
\eq{This section compares numerically the aforementioned algorithms for approximating  maps between infinite-dimensional function spaces.} Four test problems are considered:
\begin{enumerate}
    \item Navier-Stokes equation: the map between the forcing and the vorticity field at a later time is learned.
    \item Helmholtz equation:  the map between the (inhomogeneous) wavespeed field and the disturbance field (solution) is learned.
    \item Structural mechanics equation: the map between an applied boundary load and the 
    interior von Mises stress field is learned.
    \item Advection equation: the solution operator from initial condition to solution at a later time is learned.
\end{enumerate}
In our first two tests, the output functions are smooth, while the outputs of the latter two tests have discontinuities in the output space or its gradients.

This section is organized as follows: \Cref{ssec:NN-arch} describes our implementations of the neural networks from~\Cref{sec:NN}. \Cref{ssec:results} details, in \Cref{sssec:NS,sssec:Helmholtz,sssec:Solid,sssec:Advection}, each of the above four test problems; in the same subsection results are presented for each of the test problems, commenting
qualitatively on the results for each test problem in its subsection. Then, \Cref{ssec: discussion} discusses our quantitative results for all test problems,
highlighting the cost-accuracy trade-off for learning operators by neural networks.

\subsection{Neural Network Architectures}\label{ssec:NN-arch}
\begin{figure}[ht]
    \centering
    \includegraphics[width=\textwidth]{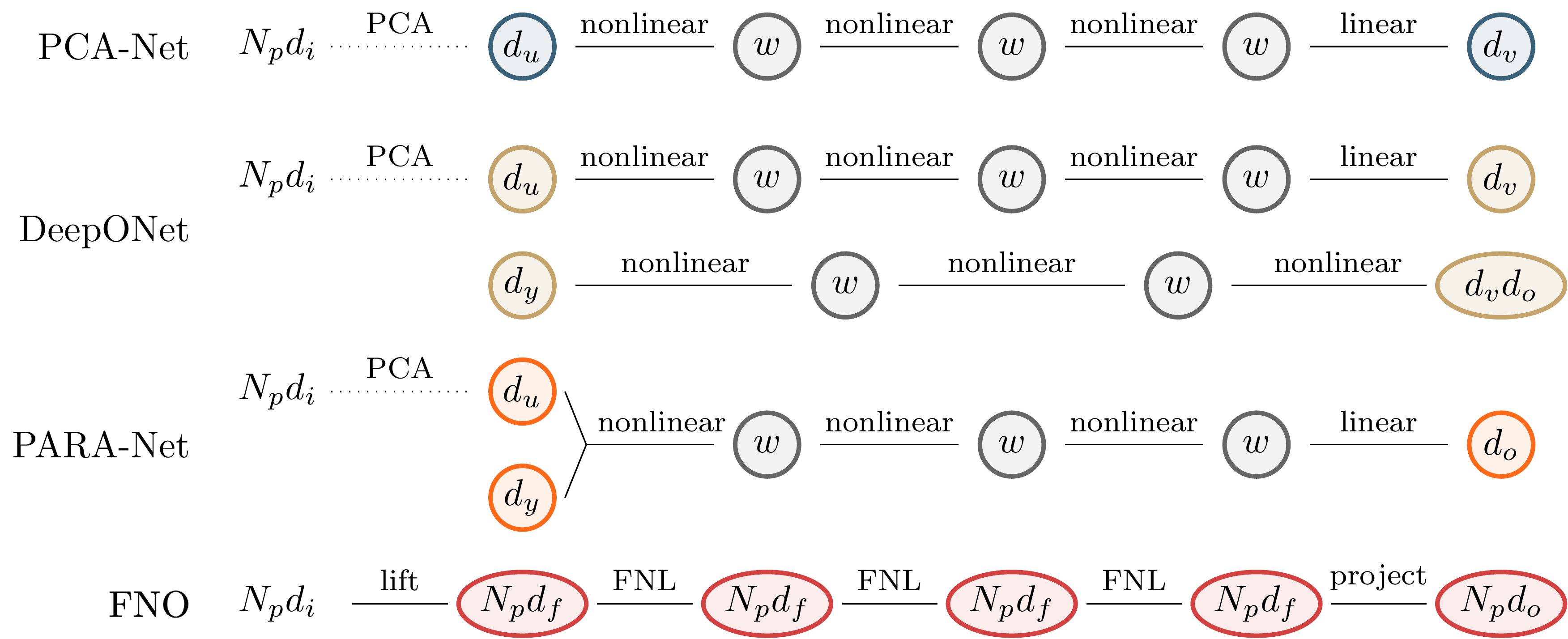}
    \caption{Schematic of neural network architectures in numerical experiments. 
    Circles represent layers; the width of each layer is given in the circle. Edges represent transformations between layers; the type of transformation between each layer is noted above each edge. Nonlinear and linear transformations are standard fully-connected layers; the lift and project layers are defined in \eqref{eq: FNO lift downsample defs}; the Fourier Neural Layer (FNL) is defined in \eqref{eq: FNL definition}. \label{fig:nn-schematic}}
\end{figure}

\Cref{fig:nn-schematic} summarizes the neural network architectures considered in our numerical experiments. In particular, the neural networks used in PCA-Net, PARA-Net, and the DeepONet branch and trunk networks have shared internal structure: they each use fully connected neural networks with three internal layers of constant fixed width $w$ between the input and output layers, and ReLU functions are employed. The output layer for each of these networks is linear (no nonlinear activation function). The FNO network has three internal Fourier Neural Layers (defined in~\eqref{eq: FNL definition}) in between the initial lifting layer and the final projection layer~\eqref{eq: FNO lift downsample defs}\eq{, and uses Gaussian Error Linear Unit (GELU) activation functions.}
PCA-Net, DeepONet and PARA-Net are initialized based on the method described in~\cite{he2015delving} and FNO is initialized following ~\cite{glorot2010understanding,li2020fourier}. 

To compare the online evaluation costs of the neural networks, we provide a cost analysis in terms of the requisite floating point operations (FLOPs) in 
\Cref{app: eval complex}; the resulting costs are tabulated in~\Cref{tab: eval complex}. DeepONet and PCA-Net have the same evaluation cost because the DeepONet branch and PCA-Net have the same architecture in our implementations, and the DeepONet trunk defines basis functions which can be precomputed after the network is trained, before evaluating the trained network on new data.

\begin{table}[h]
\centering

\begin{tabularx}{\textwidth}{l|X|c}
     Architecture & Evaluation cost & Cost scaling\\\hline
     PCA-Net & $d_u(2N_pd_i-1) + 2d_uw + 4w^2 + 2d_vw + 3w + (2d_v - 1)N_pd_o$  & $\bigO(N_p + w^2)$\\
     DeepONet & $d_u(2N_pd_i-1) + 2d_uw + 4w^2 + 2d_vw + 3w + (2d_v - 1)N_pd_o$ & $\bigO(N_p + w^2)$\\
     PARA-Net & $d_u(2N_pd_i-1) + [2(d_u+d_y)w + 4w^2 + 2wd_o + 3w]N_p$ & $\bigO(N_p w^2)$\\
     FNO & $2N_pd_f(d_i+d_o) + 3(10d_fN_p  log(N_p) + k_{max}(2d_f^2 - d_f) + 2d_f^2N_p)$ & $\bigO(d_fN_plog(N_p) + N_p d_f^2)$
\end{tabularx}

\caption{Evaluation cost of the four neural network architectures considered in this work.}\label{tab: eval complex}
\end{table}

\eq{We will} study the cost and accuracy of the neural network approximations as the amount of data and the \eq{network size} are varied. For PCA-Net, DeepONet, and PARA-Net, the widths tested are $w = \{16, 64, 128, 256, 512\}$, and the numbers of channels tested in FNO are $d_f = \{2,4, 8, 16, 32\}$. Evaluation cost is measured in terms of FLOPs (\Cref{tab: eval complex}).  
To quantify the expressive power of each of these networks, we provide a parameter complexity analysis of \eq{the} networks in \Cref{tab: param complex}. \eq{Details can be found in \Cref{app: param complex}.}

\begin{table}[ht]
\centering
\begin{tabular}{l|c}
     Architecture & Parameter complexity \\\hline
     PCA-Net & $2w^2 + w(d_u + d_v) + 3w + d_v$ \\
     DeepONet & $4w^2 + w(d_u+d_v+d_y+d_vd_o) + 6w + d_v + d_vd_o$ \\
     PARA-Net & $2w^2 + w(d_o+d_u+d_y) + 3w + d_o$\\
     FNO & $d_fd_i + d_f + d_fd_o + d_o + 3(d_f^2 + d_f^2k_{max})$
\end{tabular}
\caption{Parameter complexity of the four neural network architectures considered in this work in terms of input and output space dimensions as well as network size parameter $w$ or $d_f$.}\label{tab: param complex}
\end{table}

\subsection{Test Problems and Qualitative Results}\label{ssec:results}
This section introduces the four test problems and presents the results of our numerical comparisons of the four network architectures for each test problem. \eq{Sections \ref{sssec:NS}, \ref{sssec:Helmholtz}, \ref{sssec:Solid} and \ref{sssec:Advection}, respectively, introduce the Navier-Stokes, Helmholtz, solid mechanics, and one-dimensional advection test problems, respectively.}
For each test problem, we show comparisons of the true vs.\ neural network predicted fields at inputs that result in median and worst-case test errors and comment qualitatively on what these comparisons show about the neural network performance for each problem. Detailed quantitative discussion of the cost-accuracy trade-offs of the networks is deferred to \Cref{ssec: discussion}.

\subsubsection{Navier-Stokes Equation}
\label{sssec:NS}

\begin{figure}[t]
     \centering
     \includegraphics{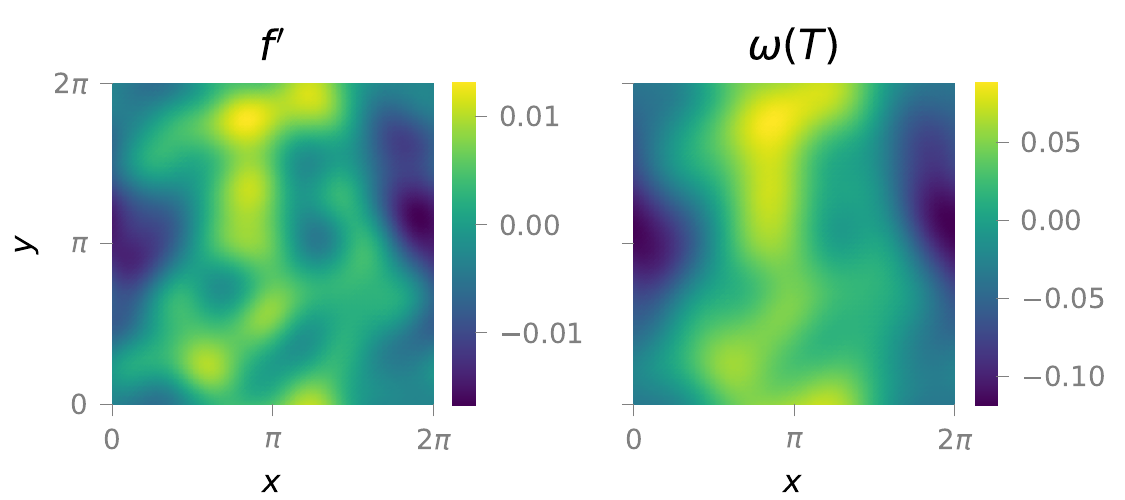}\\
     \caption{Navier-Stokes problem: sample input and output functions (left and right, respectively).}
     \label{fig:NS-map}
\end{figure}
\paragraph{Formulation}
We consider the vorticity-stream~($\omega-\psi$) formulation of the incompressible Navier-Stokes equations on a two-dimensional periodic domain, $D=D_u = D_v = [0,2\pi]^2$:
\begin{equation}
\begin{split}
    &\frac{\partial \omega}{\partial t} + (v\cdot\nabla)\omega - \nu\Delta\omega = f', \\
    &\omega = -\Delta\psi \qquad \int_{D}\psi = 0,\\
    &v = \Big(\frac{\partial \psi}{\partial x_2}, -\frac{\partial \psi}{\partial x_1}\Big).
\end{split}\label{eq: Navier Stokes vorticity stream}
\end{equation}
We are interested in the map from the forcing $f'$ to the vorticity field $\omega$ at time $t=T$.
The forcing $f'$ is assumed to be a centred Gaussian with covariance
$$\mathsf{C} = (-\Delta + \tau^2 )^{-d};$$
here $-\Delta$ denotes the Laplacian on  $D$ subject to periodic boundary conditions on the space of spatial-mean zero functions, $\tau = 3$ denotes the inverse length scale of the random field and \eq{$d = 4$ determines its regularity; the choice of $d$ then leads to up to $3$} fractional derivatives for samples from this measure.
The initial condition $\omega(0)$ is fixed and generated from the same distribution.

Our numerical experiments are conducted in the case where $\nu=0.025$, similar to the setup in
\cite{farazmand2017variational}\cite[Chapter~2.2]{majda2006nonlinear}, and we use final time $T=10.$ For the value of $\nu$ used here, the solution at time $10$ has decorrelated from the initial condition; it inherits the spatial pattern of the forcing $f'$ but has larger amplitude, and smoother small scale features.
\Cref{eq: Navier Stokes vorticity stream} is solved using a pseudo-spectral method on a $64 \times 64$ grid. To eliminate aliasing error, the Orszag 2/3-Rule~\cite{orszag1972numerical} is applied and, therefore there are $42^2$ Fourier modes (padding with zeros). Time-integration is performed using the Crank–Nicolson method with $\Delta t=10^{-3}$. 
See \Cref{fig:NS-map} for a visualization of sample input $f'$ and resultant output $\omega(\cdot,T)$ fields.

\begin{figure}[h]
     \centering
     \includegraphics[width=0.49\textwidth]{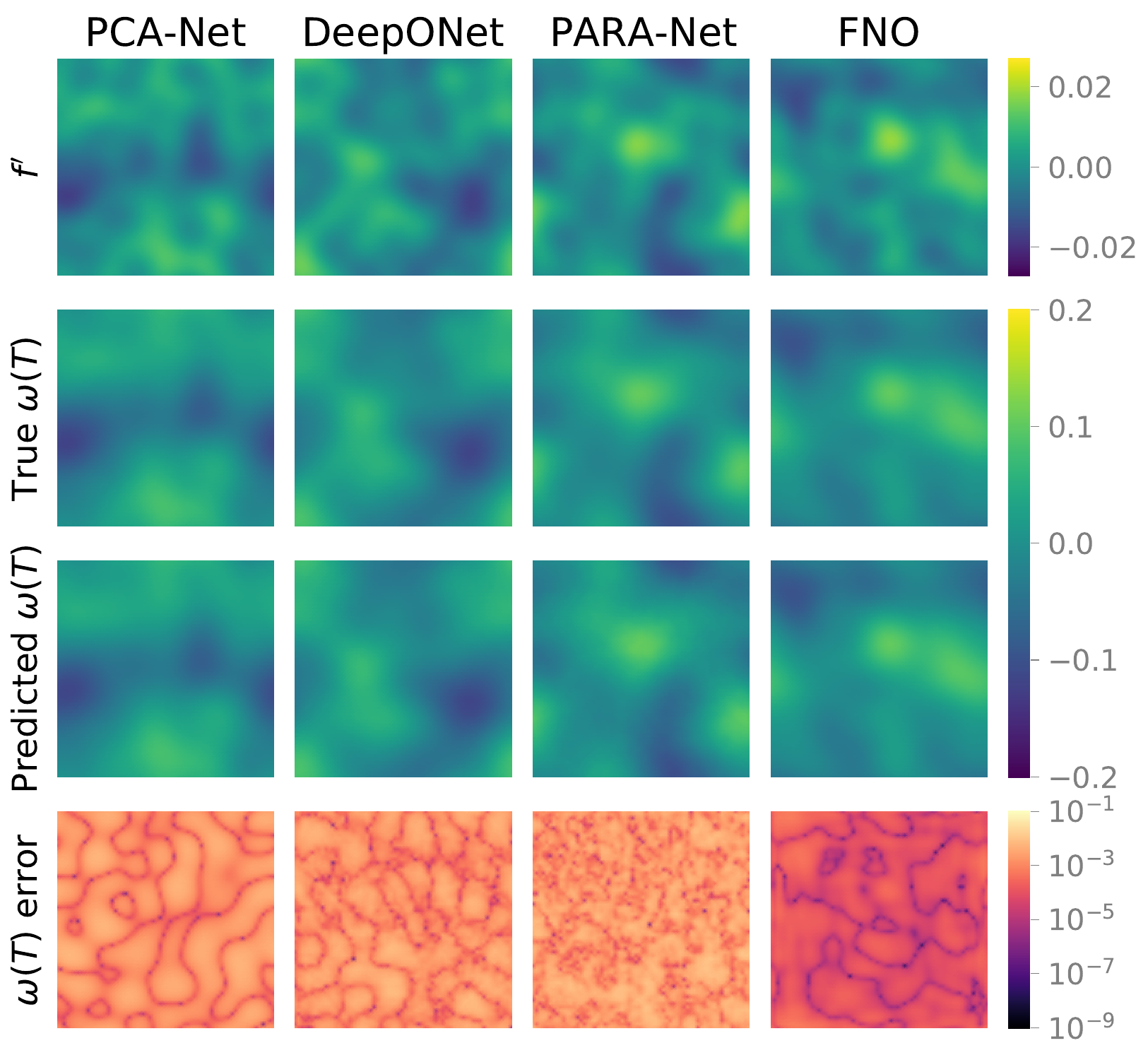}
     \includegraphics[width=0.49\textwidth]{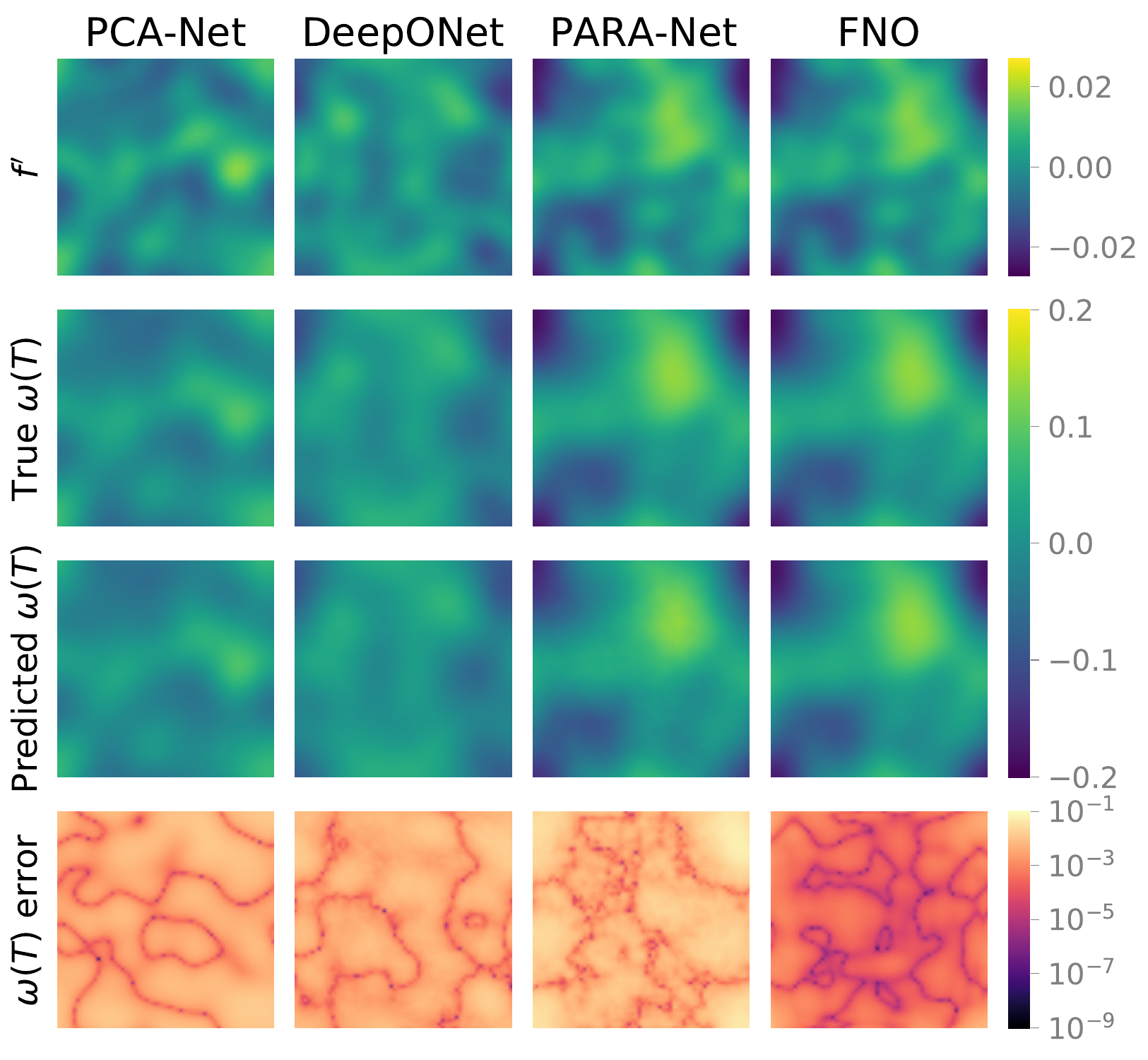}
     \caption{Navier-Stokes test problem: learned model vorticity predictions for inputs resulting in \textbf{median (left)} and \textbf{largest (right)} test errors for networks of size $w = 128$ / $d_f = 16$ trained on $N = 10000$ data. }
     \label{fig:NS-median}
\end{figure}

\paragraph{Results}
\Cref{fig:NS-median} shows the input, true output, neural network-predicted output, and output errors for the inputs resulting in the median and largest test errors (left and right) for each network architecture. 
The output $\omega(T)$ for this problem is well-correlated with the input $f'$, and all neural networks succeed in predicting the main features of the vorticity field. Note that the vorticity fields predicted by PARA-Net are grainier than those of other networks, which is reflected by the smaller length scale of the PARA-Net error fields. This graininess is due to the pointwise prediction of the network. For this problem, FNO errors are significantly lower than those of other methods, reflecting that
the problem specification is particularly well-adapted to a spectral representation.

\clearpage

\subsubsection{Helmholtz Equation}\label{sssec:Helmholtz}

\begin{figure}[t]
    \centering
    \resizebox{4.1cm}{!}{
    \begin{tikzpicture}
         \draw (0,0) -- (0,1.3) ;
         \draw (0,1.7) -- (0,3) ;
         
         \draw (0,3) -- (1.2,3) ;
         \draw (1.8,3) -- (3,3) ;
         
         \draw (3,3) -- (3,1.7) ;
         \draw (3,1.3) -- (3,0) ;
         
         \draw (0,0) -- (1.2,0) ;
         \draw (1.8,0) -- (3,0) ;

         \node at (1.5,1.5) {$\Omega$};
         \node at (1.5,0) {\eight $\partial\Omega_1$};
         \node at (3,1.5) {\eight $\partial\Omega_2$};
         \node at (1.5,3) {\eight $\partial\Omega_3$};
         \node at (0,1.5) {\eight $\partial\Omega_4$};
         \node at (1.5,-0.5) [color=white]{$x$};
    \end{tikzpicture}
    }
    \includegraphics[height=0.3\textwidth]{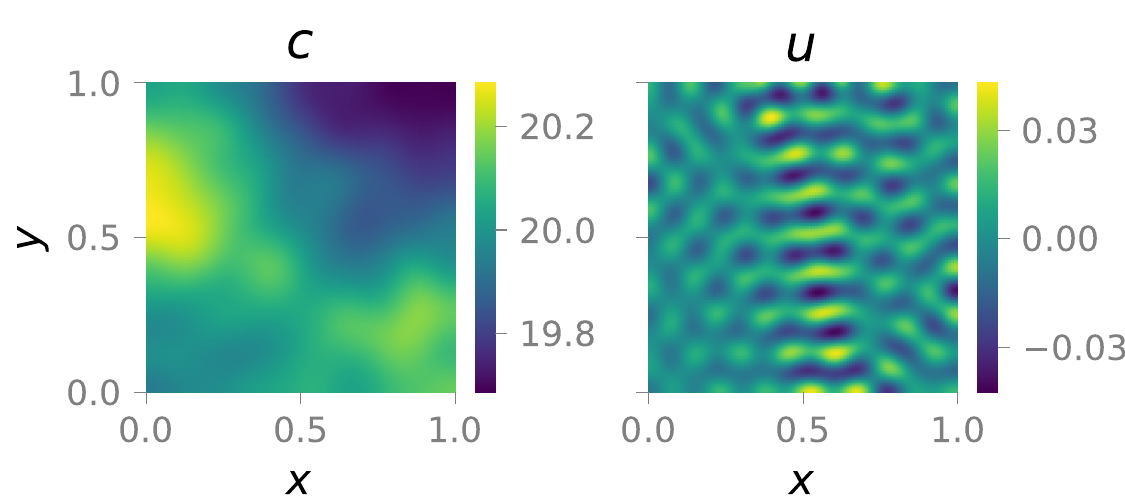}
    \caption{Helmholtz test problem: \textbf{Left:} schematic of unit domain with labeled boundaries. \textbf{Center: } Sample input wave speed field. \textbf{Right: }Sample output disturbance field.
    \label{fig: Helmholtz domain}}
\end{figure}
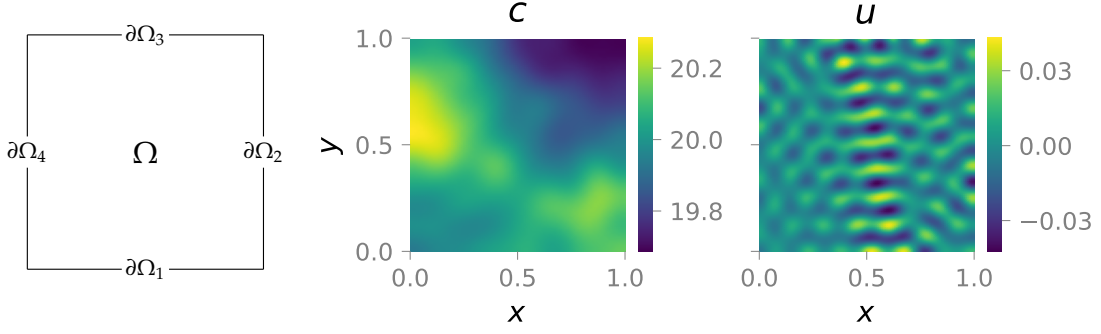

\paragraph{Formulation}
We consider the Helmholtz equation on the domain $D=D_u = D_v = [0,1]^2$
shown in \Cref{fig: Helmholtz domain}. Given frequency $\omega = 10^3$
and wavespeed field $c: \Omega \to \R$, the excitation field 
$u: \Omega \to \R$ solves equation
\begin{subequations}
  \label{eq: Helmholtz BVP}
  \begin{align}
      \Bigl(-\Delta - \frac{\omega^2}{c^{2}(x)}\Bigr) u &= 0  & &\text{ in } \Omega,\\
      \frac{\partial u}{\partial n} &= 0  & &\text{ on } \partial\Omega_1,\ \partial\Omega_2,\ \partial\Omega_4 \\
       \frac{\partial u}{\partial n} &= u_{N}  & &\text{ on } \partial\Omega_3.
  \end{align}
\end{subequations}
Note that the Neumann boundary condition imposed on $\partial \Omega$ is non-zero only on the top edge. Throughout the experiments
presented here $u_{N}$ is fixed at $1_{\{0.35 \leq x \leq 0.65\}}$. The wavespeed field $c(x)$ is assumed to be 
$$c(x) = 20 +  \textrm{tanh} (\tilde{c}(x)),$$ 
where $ \tilde{c}$ is a centred Gaussian 
$$\tilde{c} \sim \N(0, \mathsf{C}) \quad \textrm{and} \quad \mathsf{C} = (-\Delta + \tau^2 )^{-d};$$
here $-\Delta$ denotes the Laplacian on $D_u$ subject to homogeneous Neumann boundary conditions on the space of spatial-mean zero functions, and  we choose $d=2$ and $\tau = 3$,
the choice of $d$ then leads to up to $1$ fractional derivative for samples from this measure.
We are interested in the map from the wavespeed field $c$ to the solution $u$. Equation~\eqref{eq: Helmholtz BVP} is solved using a finite element method on a $100\times100$ grid. See \Cref{fig: Helmholtz domain} for a visualization of sample input (c) and resultant output $u$ fields.

\begin{figure}[h]
     \centering
    \includegraphics[width = 0.49\textwidth]{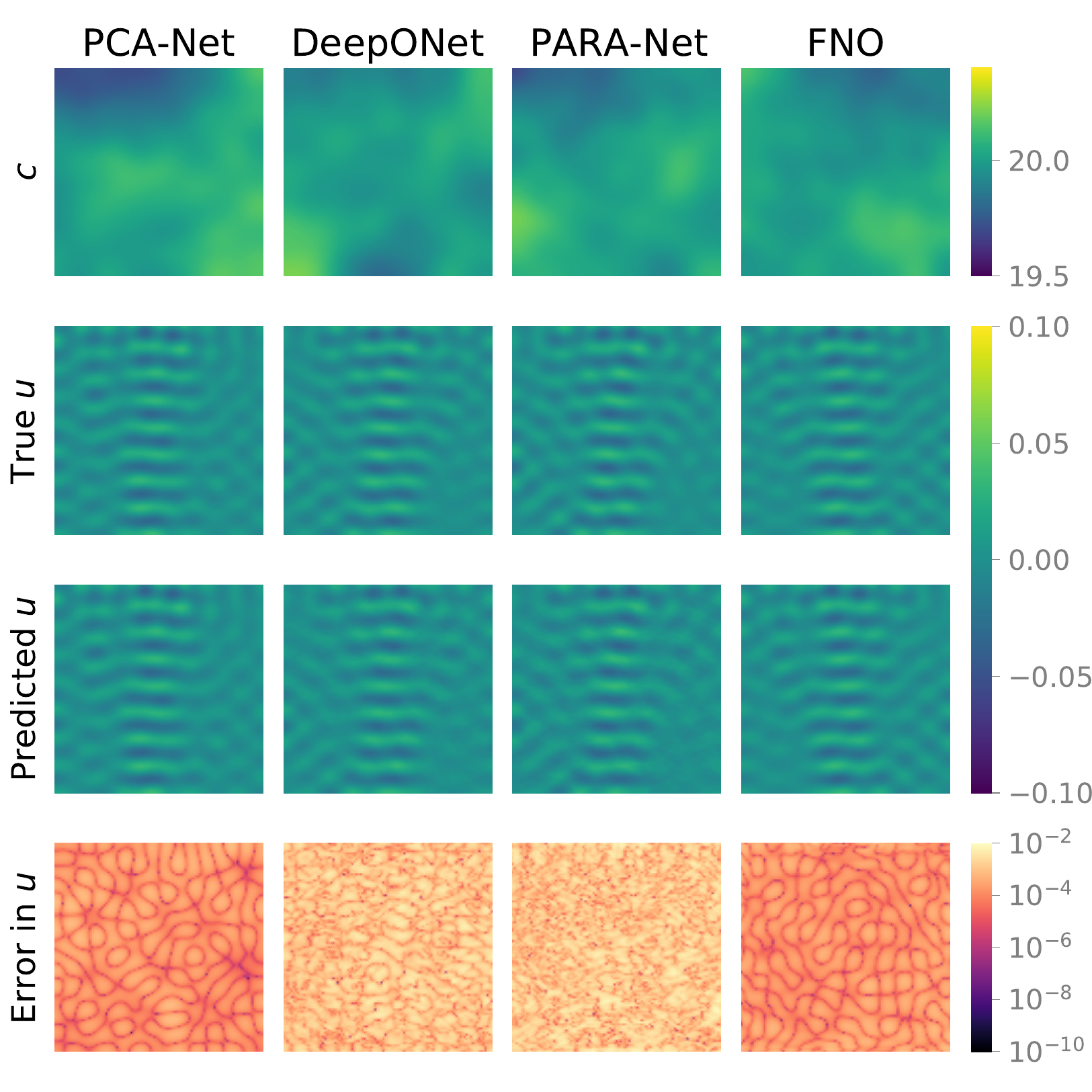}\hfill
    \includegraphics[width=0.49\textwidth]{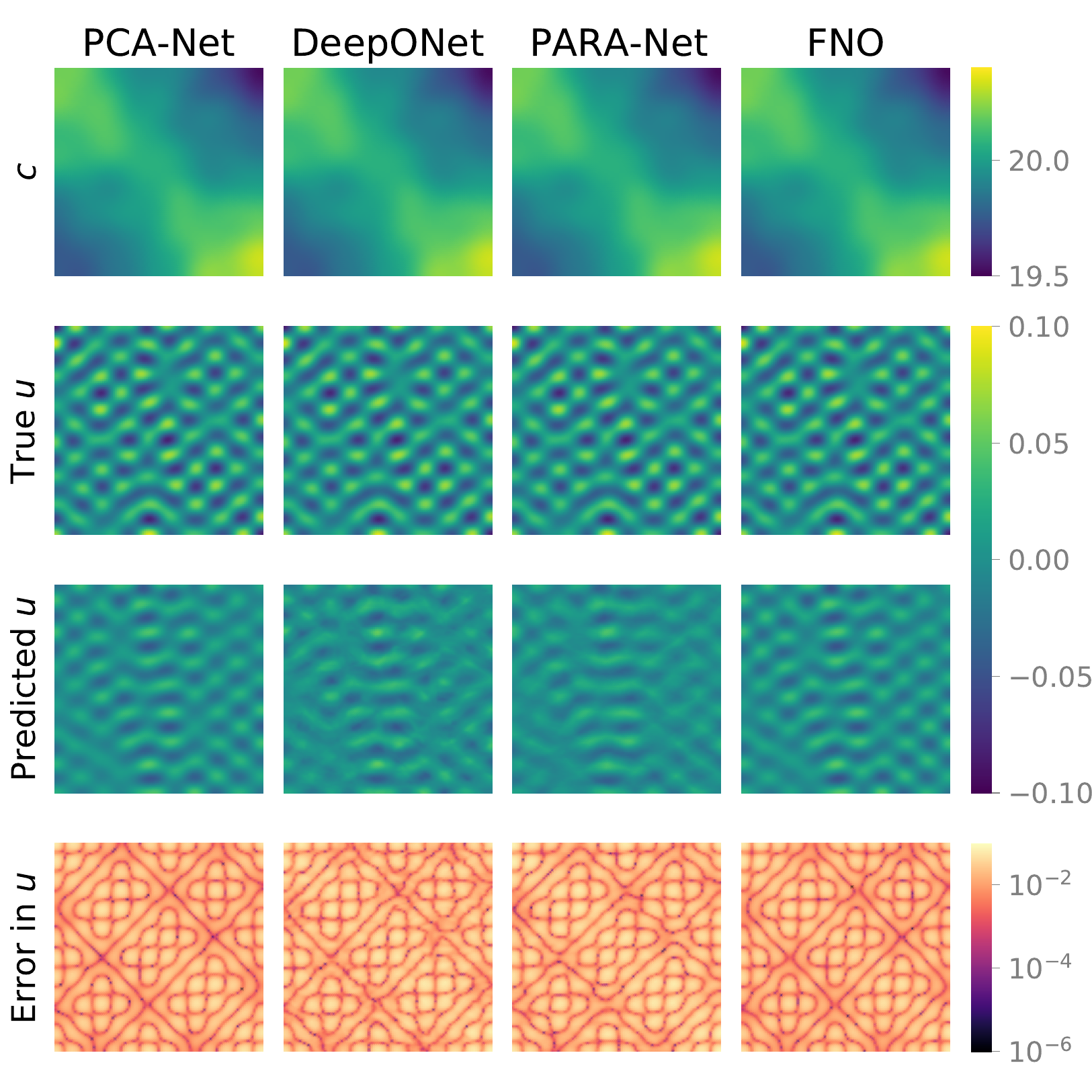}
     \caption{Helmholtz test problem: learned model predictions for inputs resulting in \textbf{median (left)} and \textbf{largest (right)} test errors for networks of size $w = 128$ / $d_f = 16$ trained on $N = 10000$ data.}
     \label{fig:Helmholtz-median}
\end{figure}

\paragraph{Results}
\Cref{fig:Helmholtz-median} shows the input, true output, neural network-predicted output, and output errors for the inputs resulting in the median and largest test errors (left and right panels) for each network architecture. 
For the median error cases, all neural networks tested yield disturbance predictions that match the true disturbance field in the `eyeball norm'. However, differences in their behavior are revealed by the structure of the error fields: the error fields of PCA-Net and FNO have similar structures, length scales, and error magnitudes, indicating that their output spaces are similar. In contrast DeepONet and PARA-Net exhibit error fields with smaller length scales and different structures.

\begin{figure}[h]
    \centering
    \includegraphics[width=0.95\textwidth]{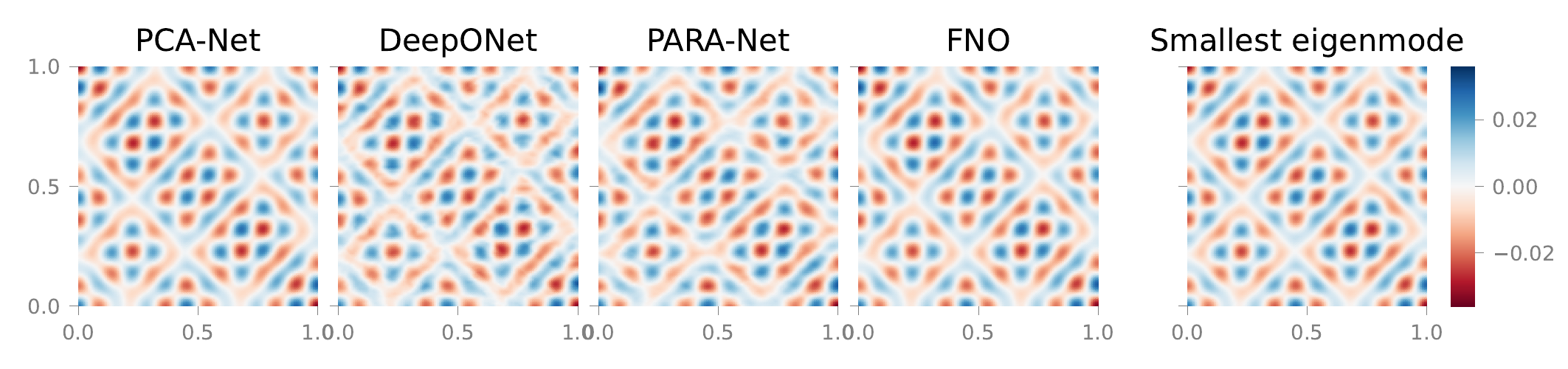}
    \caption{Helmholtz test problem largest test error case: normalized differences between learned model predictions and true solution (left four panels) and comparison to normalized smallest eigenmode of the Helmholtz operator for this test case. }
    \label{fig:helmholtz-eigen}
\end{figure}

For the worst-case error cases, \eq{we note that all four networks achieve the worst-case error on the same test input $c$. This is because for this particular $c$, the frequency $\omega$ is close to an eigenfrequency of the Helmholtz operator, for which the $c\mapsto u$ map would not be uniquely defined. In~\Cref{fig:helmholtz-eigen}, we plot the normalized differences between the network predictions and the true solution, as well as the normalized eigenmode corresponding to the smallest eigenvalue of the Helmholtz operator for this input (the last row of the right panel of \Cref{fig:Helmholtz-median} corresponds to the absolute value of the left four panels of \Cref{fig:helmholtz-eigen} plotted on a log scale). Note that these differences are close to scalar multiples of the eigenmode. A sample, $c$, drawn from the distribution for testing clearly can result in an eigenfrequency close the frequency already chosen for the simulations. The solution becomes large near eigenfrequencies, as is quantified, for example, in \cite[Section~2.1]{beretta2016inverse}; hence, the accuracy of any prediction deteriorates accordingly independent of the architecture. This effect is visible in \Cref{fig:Helmholtz-median}: note that the left and right panels share a color scale and that the true $u$ for the median error has more washed out colors than the true $u$ for the largest error.}

\pagebreak
\subsubsection{Structural Mechanics Equation}\label{sssec:Solid}
\begin{figure}[t]
    \centering
    \resizebox{3.9cm}{!}{
    \begin{tikzpicture}
         \draw (0,0) -- (0,3) ;
         \draw (0,3) -- (3,3) ;
         \draw (3,3) -- (3,0) ;
         \draw (0,0) -- (3,0) ;
         \fill [pattern = north east lines] (0.0,-0.3) rectangle (3.0, 0.0);
         \draw (1.5,1.5) circle (0.9) node [below] {};
         \node at (2.0,2.5) {$\Omega$};
          \draw[thick,->] (0, 3) -- (0, 3.5);
        \draw[thick,->] (0.5, 3) -- (0.5, 3.5);
        \draw[thick,->] (1, 3) -- (1, 3.5) node[anchor=west]{$\bar{t}$};
        \draw[thick,->] (1.5, 3) -- (1.5, 3.5);
        \draw[thick,->] (2.0, 3) -- (2.0, 3.5);
        \draw[thick,->] (2.5, 3) -- (2.5, 3.5);
        \draw[thick,->] (3.0, 3) -- (3.0, 3.5);
    \end{tikzpicture}
    }
    \includegraphics[height=1.96in]{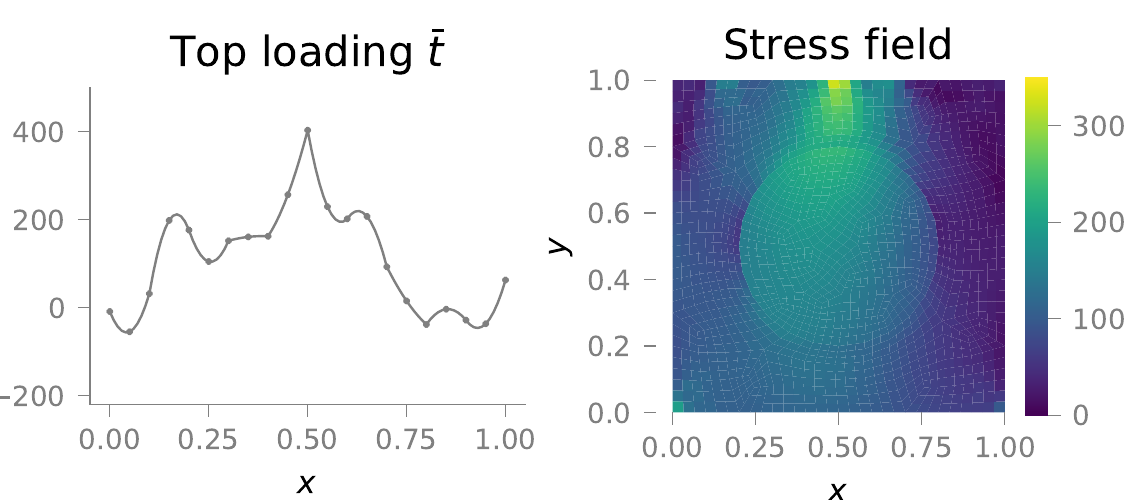}
    \caption{Solid mechanics test problem. Schematic of fiber-reinforced material and loading (\textbf{left}), sample load drawn from Gaussian distribution (\textbf{center}), and corresponding von Mises stress field (\textbf{right}).\label{fig: Solid domain}}
\end{figure}
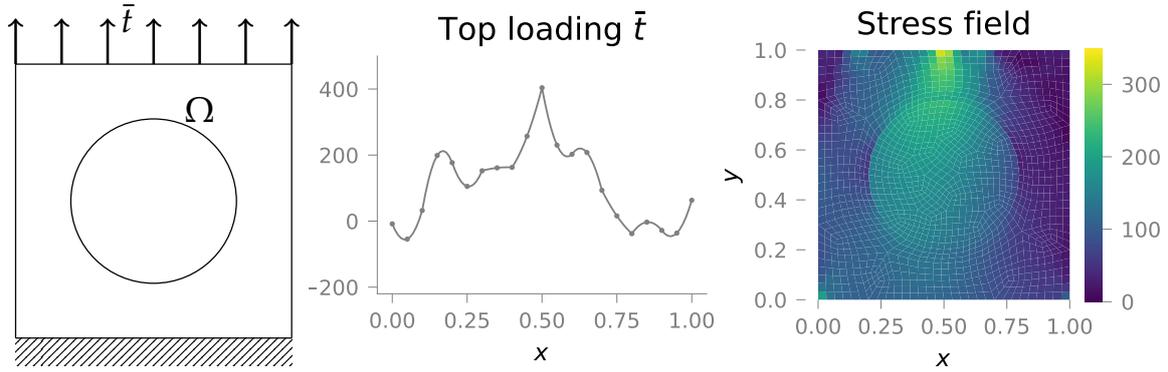

\paragraph{Formulation}
The governing equation of an elastic solid undergoing infinitesimal deformations is
\begin{equation}
  \begin{aligned}
       \nabla~\cdot \sigma & =  0 &  & \textrm{ in } \Omega,     \\
    u                & = \bar{u}   &  & \textrm{ on } \Gamma_{u}, \\
    \sigma \cdot n & = \bar{t}     &  & \textrm{ on } \Gamma_{t}, \\
  \end{aligned}
  \label{EQ:LINEAR_MOMENTUM}
\end{equation}
where $u$ is the displacement vector and $\sigma$ is the (Cauchy)~stress tensor; $\Omega$ denotes the  computational domain.
The prescribed displacement $\bar{u}$ and the surface traction $\bar{t}$ respectively, are imposed on the domain boundaries $\Gamma_{u}$ and $\Gamma_{t}$ respectively, with the outward unit normal $n$, where $\Gamma_u \cap \Gamma_t = \emptyset$ and $\Gamma_u \cup \Gamma_t = \partial \Omega$.
The unit cell problem is depicted in~\Cref{fig: Solid domain}, which is clamped on the bottom edges. Tension traction is applied on the top edge, and the distributed load is $\bar{t}$. To solve for the displacement $u$ from~\cref{EQ:LINEAR_MOMENTUM}, we also need the constitutive model, which maps the deformation gradient to the stress.
The matrix is made of incompressible Rivlin-Saunders material~\cite{pascon2019large} with density $\rho=0.8$ and energy density function parameters $C_1 = 1.863\times10^5, \quad C_1 = 9.79\times10^3$ and the cylindrical fiber at the center is made 
of linear elastic material with density $\rho=3.2$, Young's modulus $E=4\times10^6$ and Poisson ratio $\nu =0.35$.

We are interested in the map from the one-dimensional load $\bar{t}$ to the von Mises stress field $\tau_{vM}$ on the two dimensional domain $\Omega.$ The load $\bar{t}$ is \eq{drawn from a Gaussian random field} with mean $100$ and covariance $400^2 \mathsf{C}$ with 
$$\bar{t} \sim \N\big(100, 400^2 \mathsf{C}\big) \quad \textrm{and} \quad \mathsf{C} = (-\Delta + \tau^2 )^{-d}.$$
Here $-\Delta$ denotes the Laplacian on $D$ subject to homogeneous Neumann boundary conditions on the space of spatial-mean zero functions, 
$\tau = 3$ denotes the inverse length scale of the random field and $d = 1$ determines its regularity (upto $1/2$ a fractional derivative for samples from this measure). The data is generated using the \texttt{NNFEM} library~\cite{huang2020learning,xu2021learning}. The mesh consists of $189$ quadratic quadrilateral
elements and the top edge is discretized by $10$ quadratic elements and hence $21$ points. The inputs are the load on the $21$ points, and the outputs are the stress field~(see \Cref{fig: Solid domain}) on Gaussian quadrature points ($9 \times 189$). Since FNO operates on uniform grids and requires the input and output data have the same dimensions, the stress field is interpolated on a $41 \times 41$ grid via radial basis function interpolation, and the load is interpolated on a $41$ grid and extruded in the $y$ direction. \eq{We note that in generating the input and output data for this problem, the discretized input data for the load $\bar{t}$ are drawn from a Gaussian, but because the finite element solver employs quadratic elements the effective load in the solver is smoother (as depicted in \Cref{fig: Solid domain}).}

\begin{figure}[H]
     \centering
     \includegraphics[width = 0.8\textwidth]{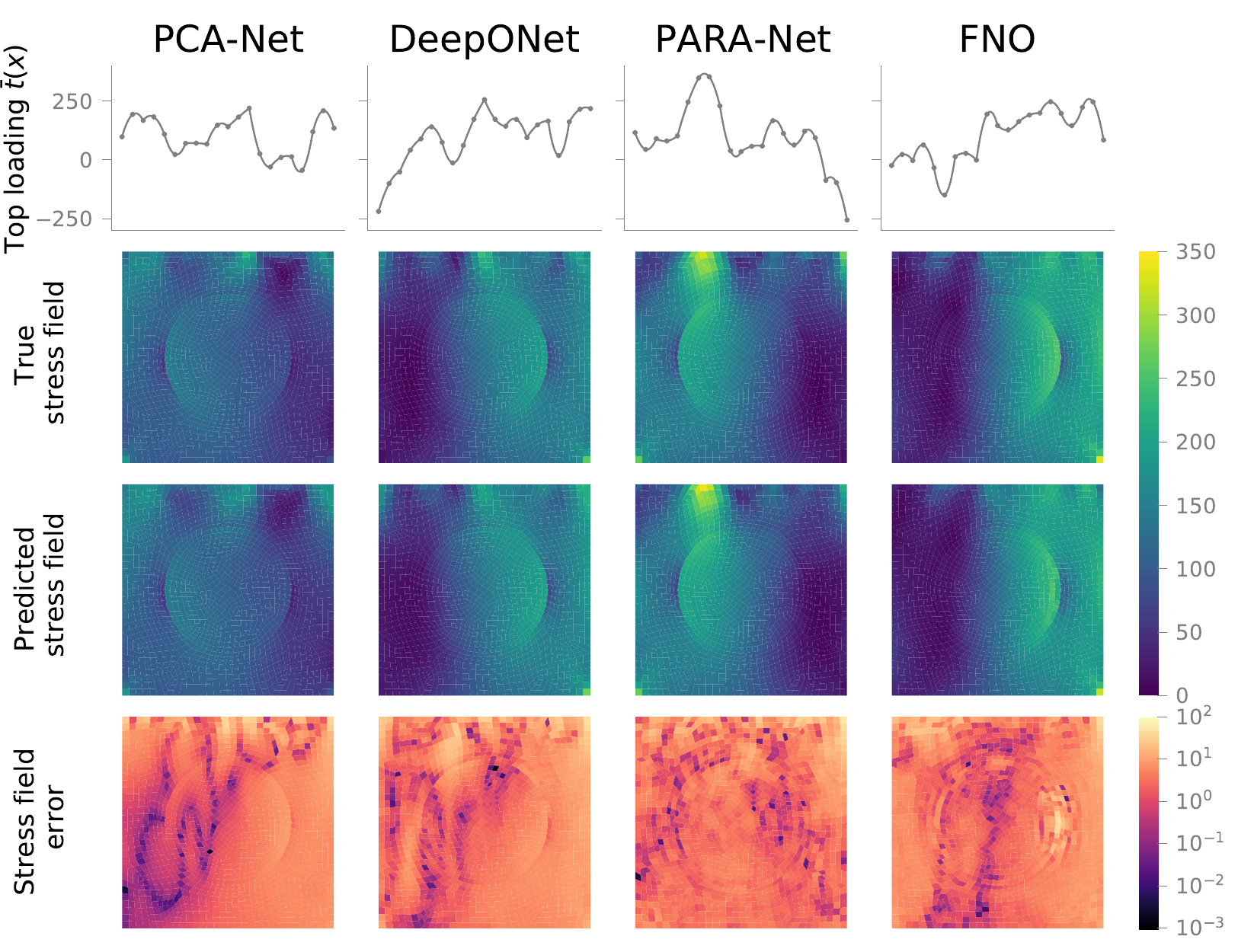}
     \caption{Structural mechanics test problem: learned model stress field predictions for inputs resulting in \textbf{median} test errors for networks of size $w = 128$ / $d_f = 16$ trained on $N = 10000$ data.}
     \label{fig:Solid-median}
\end{figure}

\begin{figure}[H]
     \centering
     \includegraphics[width = 0.8\textwidth]{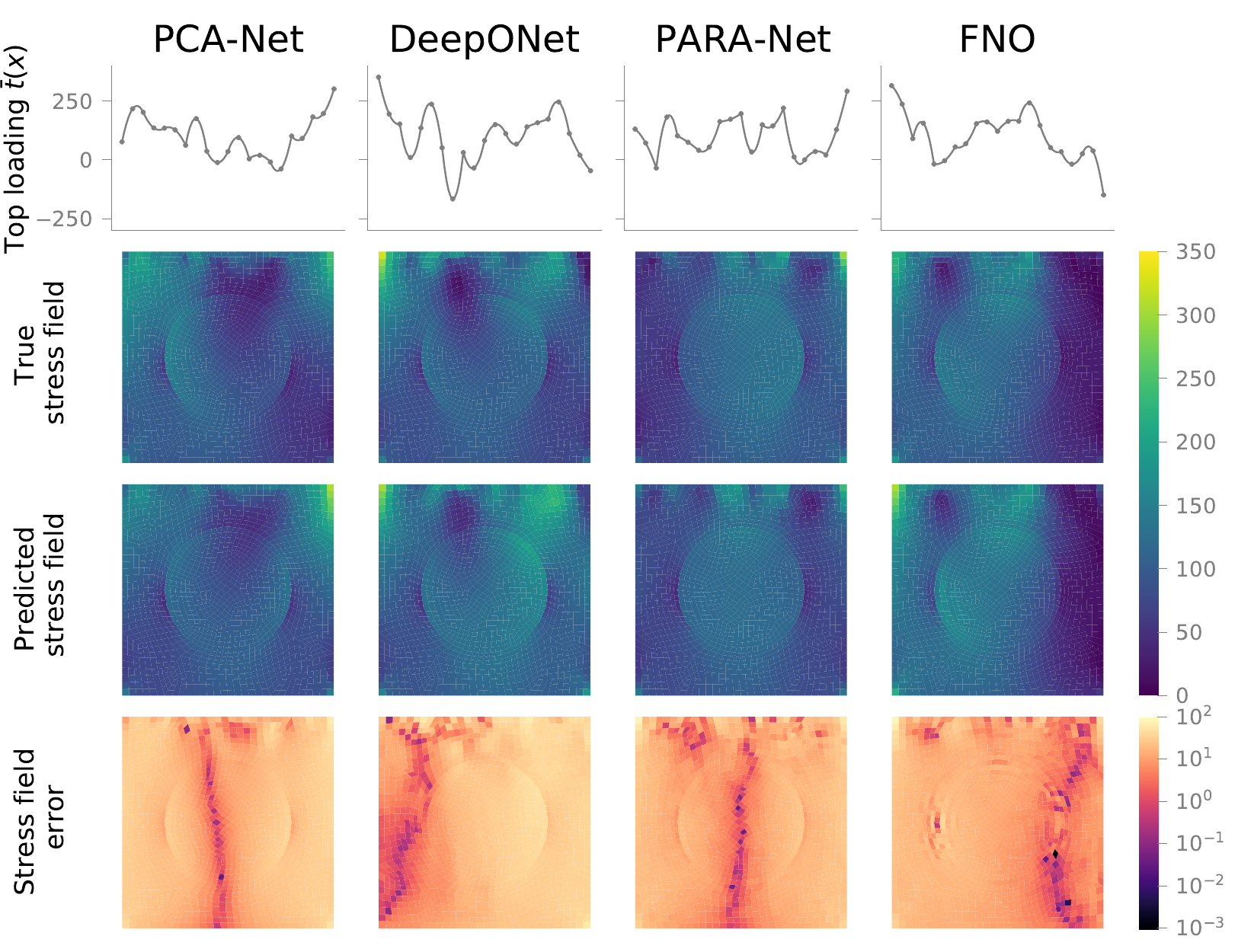}
     \caption{Structural mechanics test problem: learned model stress field predictions for inputs resulting in \textbf{largest} test errors for networks of size $w = 128$ / $d_f = 16$ trained on $N = 10000$ data.}
     \label{fig:Solid-largest}
\end{figure}

\paragraph{Results}
\Cref{fig:Solid-median,fig:Solid-largest} show the input, true output, neural network-predicted output, and output errors for inputs resulting in the median and largest test errors, respectively, for each network architecture. For this problem with discontinuous outputs, all four neural networks yield qualitatively similar predictions on both their median and worst-case error inputs, as seen in the similar color scaling of the error field plots. We call attention to the stress jump at the material interface which is generally well-captured by all four networks. We note the FNO results have oscillations at the material interface due to interpolation on the uniform grid.

\subsubsection{Advection Equation}\label{sssec:Advection}
\begin{figure}
     \centering
     \includegraphics[height=0.3\textwidth]{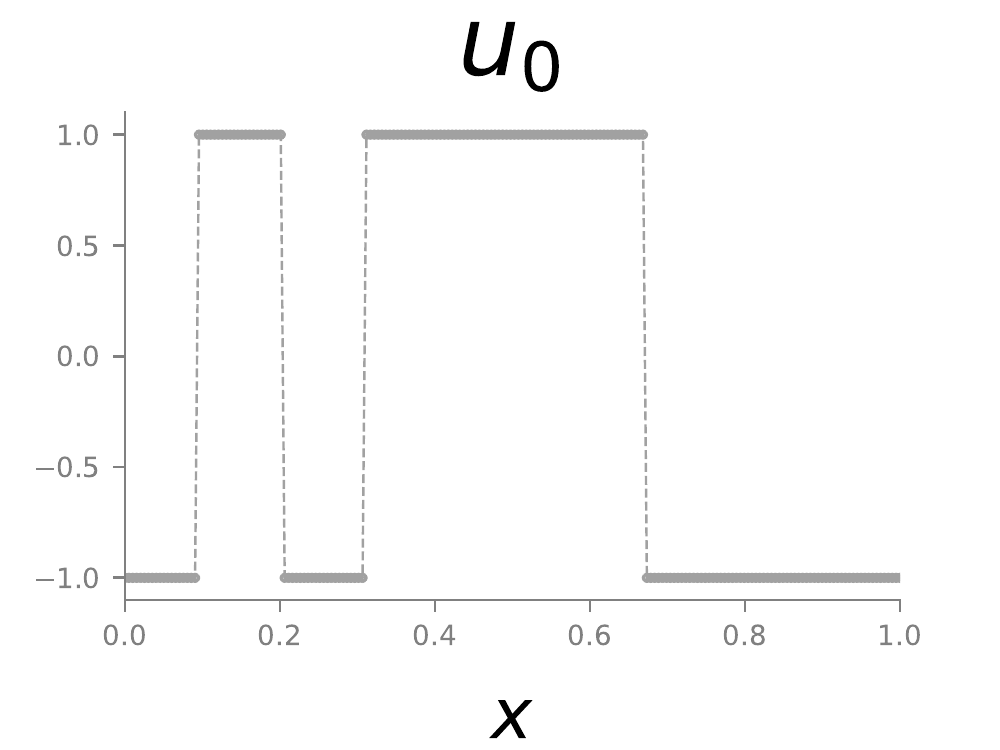}~\includegraphics[height=0.3\textwidth]{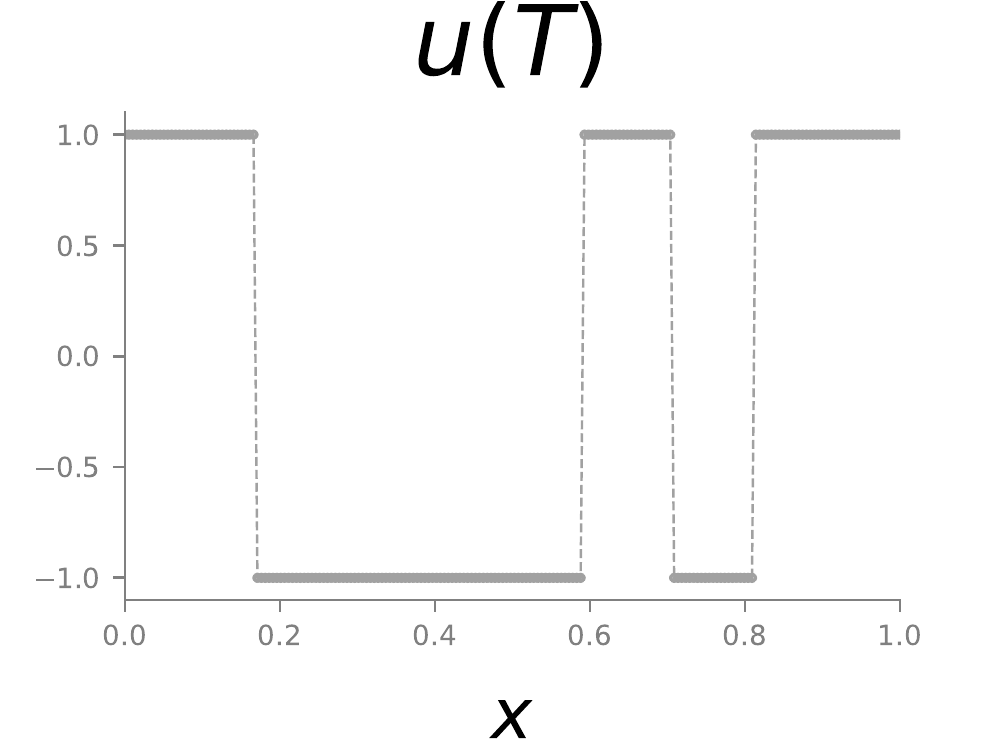}\\
     \caption{Sample input/output functions for the advection equation, left: $u_0$, right: $u(T)$.}
     \label{fig:Advection-low-in-out}
\end{figure}

\paragraph{Formulation}
The 1D advection equation in $\Omega = [0,\,1)$ is 
\begin{equation}
  \begin{aligned}
    &\frac{\partial u}{\partial t}  + c \frac{\partial u}{\partial x}  = 0  \qquad x \in \Omega, \\
    &u(0) = u_0 \\
  \end{aligned}
  \label{EQ:ADV}
\end{equation}
where $c=1$ is the constant advection speed, and periodic boundary conditions are imposed.
We are interested in the map from the initial $u_0$ to solution 
$u(\cdot,T)$ at $T = 0.5$. The initial condition $u_0$ is assumed to be 
$$u_0 = -1 + 2\1\{\widetilde{u_0} \geq 0\}$$
where $\widetilde{u_0}$ a centered Gaussian
$$\widetilde{u_0} \sim \N(0, \mathsf{C})\quad \textrm{and} \quad \mathsf{C} = (-\Delta + \tau^2 )^{-d};$$
here $-\Delta$ denotes the Laplacian on $D$ subject to periodic conditions on the space of spatial-mean zero functions, 
$\tau = 3$ denotes the inverse length scale of the random field and $d = 2$ determines 
the regularity of $\widetilde{u_0}$, which is upto $3/2$ derivatives. 
A pair of sample input and output data is depicted~\Cref{fig:Advection-low-in-out}.
\eq{Note that multiple discontinuities exist in this input sample; in general
the probability of drawing a function with a given number $p$ of discontinuities
decreases with $p$ but, in principle, draws with any number of discontinuities
are possible.}

\begin{figure}
     \centering
     \includegraphics[width=\textwidth]{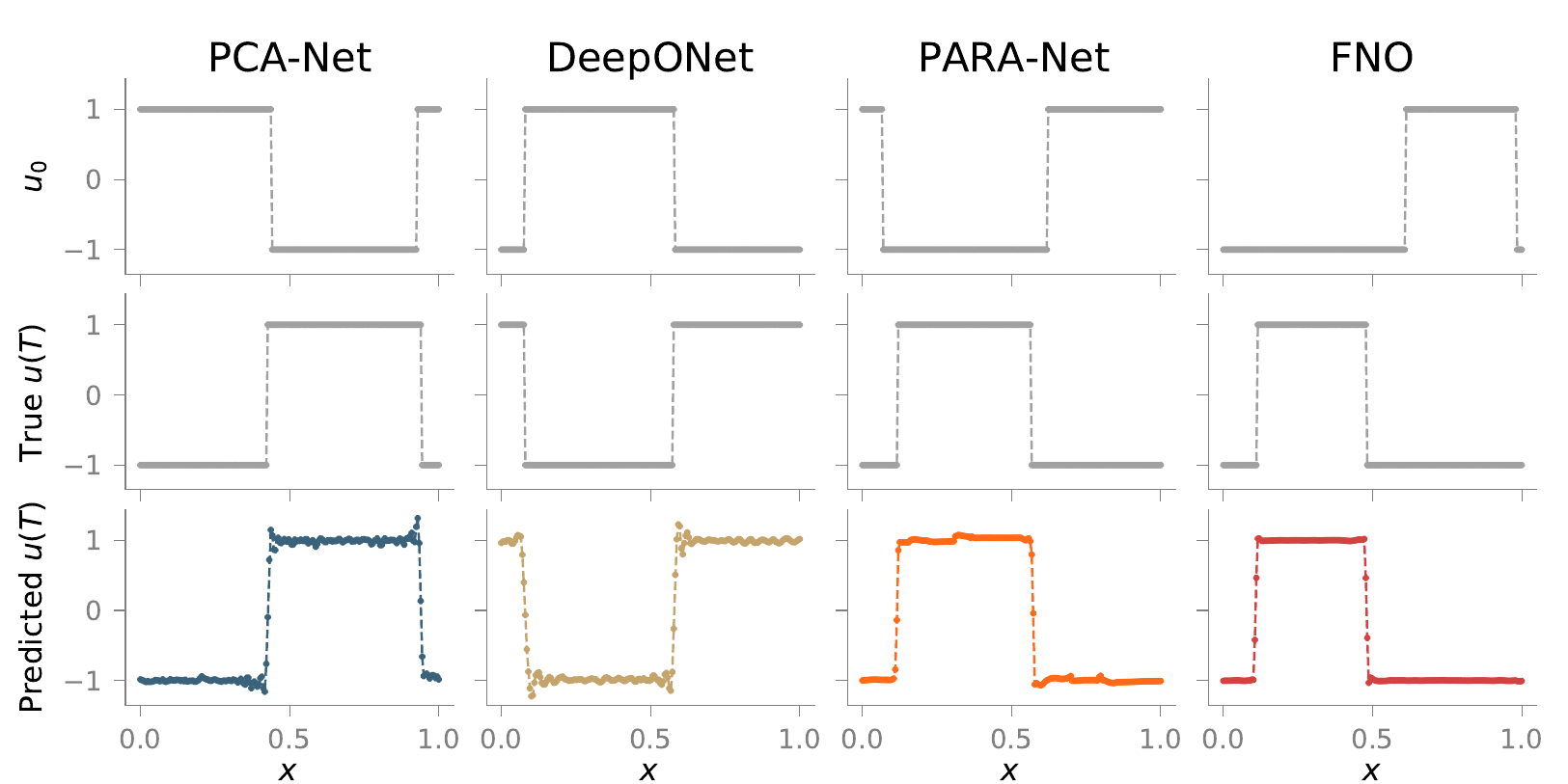}
     \includegraphics[width=\textwidth]{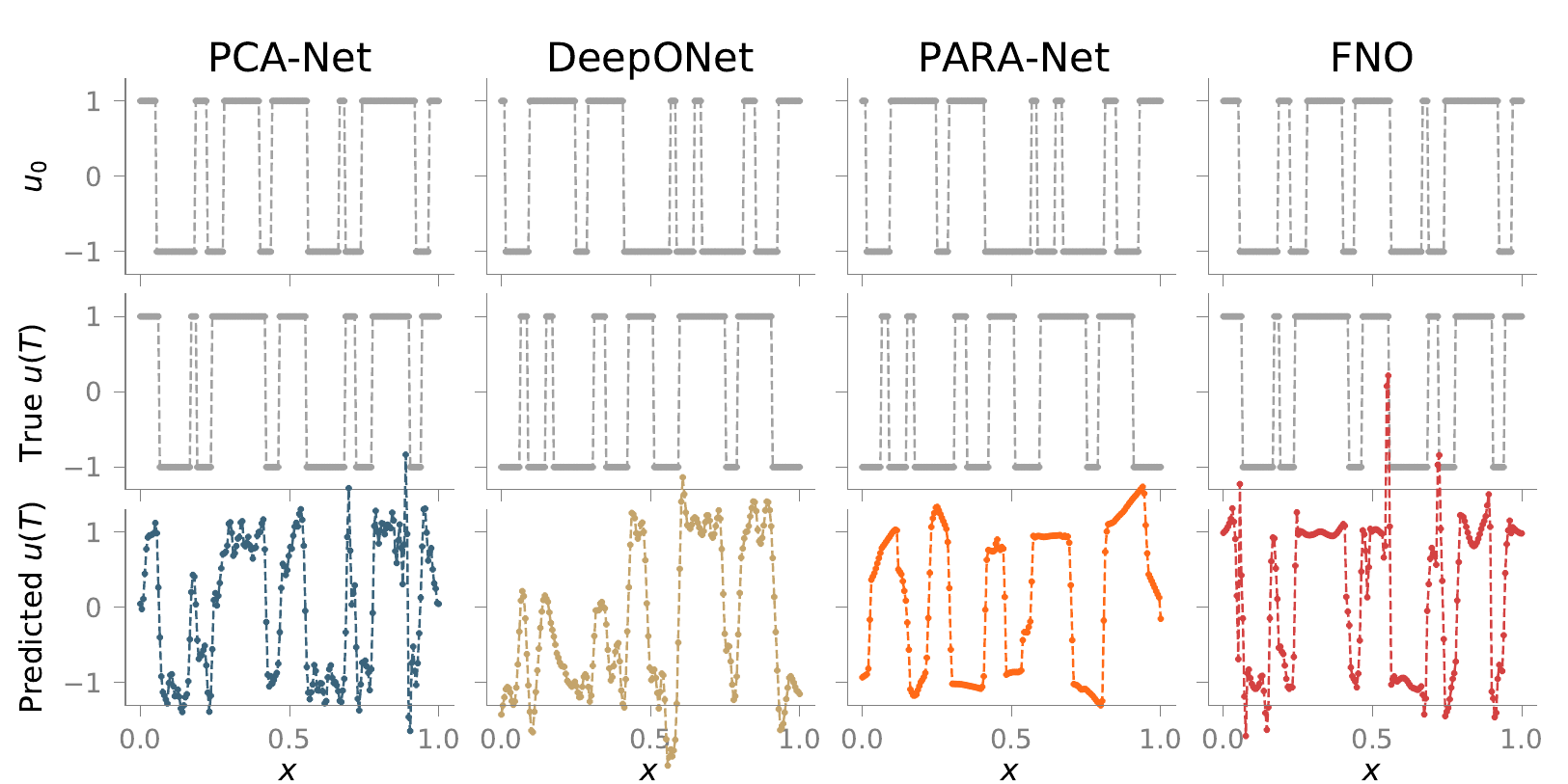}
     \caption{Advection test problem: learned model output predictions for inputs resulting in \textbf{median (top)} and \textbf{largest (bottom)} test errors for networks of size $w = 128$ / $d_f = 16$ trained on $N = 10000$ data. }
     \label{fig:Advection-low-median}
\end{figure}

\paragraph{Results}
\eq{Equation~\eqref{EQ:ADV} is solved analytically on a $200$-point uniform grid.}
\Cref{fig:Advection-low-median} shows the input, true output, and neural network-predicted output, for inputs resulting in the median and largest test errors (top and bottom panels) for each neural network architecture. We note that the median test cases for all four networks are similar; they each have just two discontinuities that are about half the domain apart. All four neural networks yield predictions which accurately reflect these discontinuities in the median case, although PCA-Net and DeepONet suffer from oscillatory Gibbs phenomena near the discontinuities. The worst-case test cases across all four network architectures also share similar characteristics: they have many discontinuities, including up-and-down jumps within a length scale of about one-tenth of the domain. \eq{In these challenging cases, the PCA-Net and DeepONet predictions suffer from Gibbs phenomena throughout the domain due to the many discontinuities. The FNO prediction reflects discontinuities at longer length scales as well and, additionally, suffers from significant overshooting near to small-length-scale discontinuities. Lastly, the PARA-Net worst-case prediction outputs a piecewise smooth solution that does not exhibit Gibbs phenomena but  does not reflect the true discontinuities of the solution. It is important to note that
all of PCA-Net, DeepONet and PARA-Net use the same input space representation,
based on projection onto a finite number of PCA modes, and this is not well-adapted
to discontinuous inputs. However DeepONet and PARA-Net have the possibility of
recovering from this, since they learn the output space representation and it is
in the output space that the error is measured.}


\subsection{Discussion of Quantitative Results}\label{ssec: discussion}
We now consider the complexity question as measured by test errors vs.\ cost for each of the neural networks. There are two main axes along which we can measure cost, training cost and online evaluation cost. For many PDE problems, the training cost is dominated by the cost of generating data by numerically solving the PDEs; we thus measure training cost in terms of available training data volume. In our experimental design, the online evaluation cost is directly related to the network size. The network size and training data volume have coupled influence on the network accuracy. The test error can thus be viewed as a surface in three-dimensional space where the two independent variables are the training volume and network size.
In \Cref{sssec: acc vs training}, we consider slices of this surface along the data volume axis: we report and discuss the error-vs-training cost curves for each network and test problem at different network sizes. Then, in \Cref{sssec: acc vs size}, we consider slices of the error in the other direction, along the network size axis, and report and discuss error-vs-size cost curves for each of the networks and test problems. \eq{\Cref{sssec: acc v cost} reports and discusses the error-vs-online evaluation cost curves. \Cref{app: out of distribution tests} provides a \eq{preliminary} study of out-of-distribution generalization error. \Cref{sssec:output} compares the output space representations learned by PCA-Net and DeepONet. Finally, in \Cref{sssec: role of opt} we comment on the influence of optimizer performance on the neural network prediction errors.}

\subsubsection{Accuracy vs.\ Training Cost}\label{sssec: acc vs training}
\begin{figure}[h]
     \centering
     \includegraphics[width=0.95\textwidth]{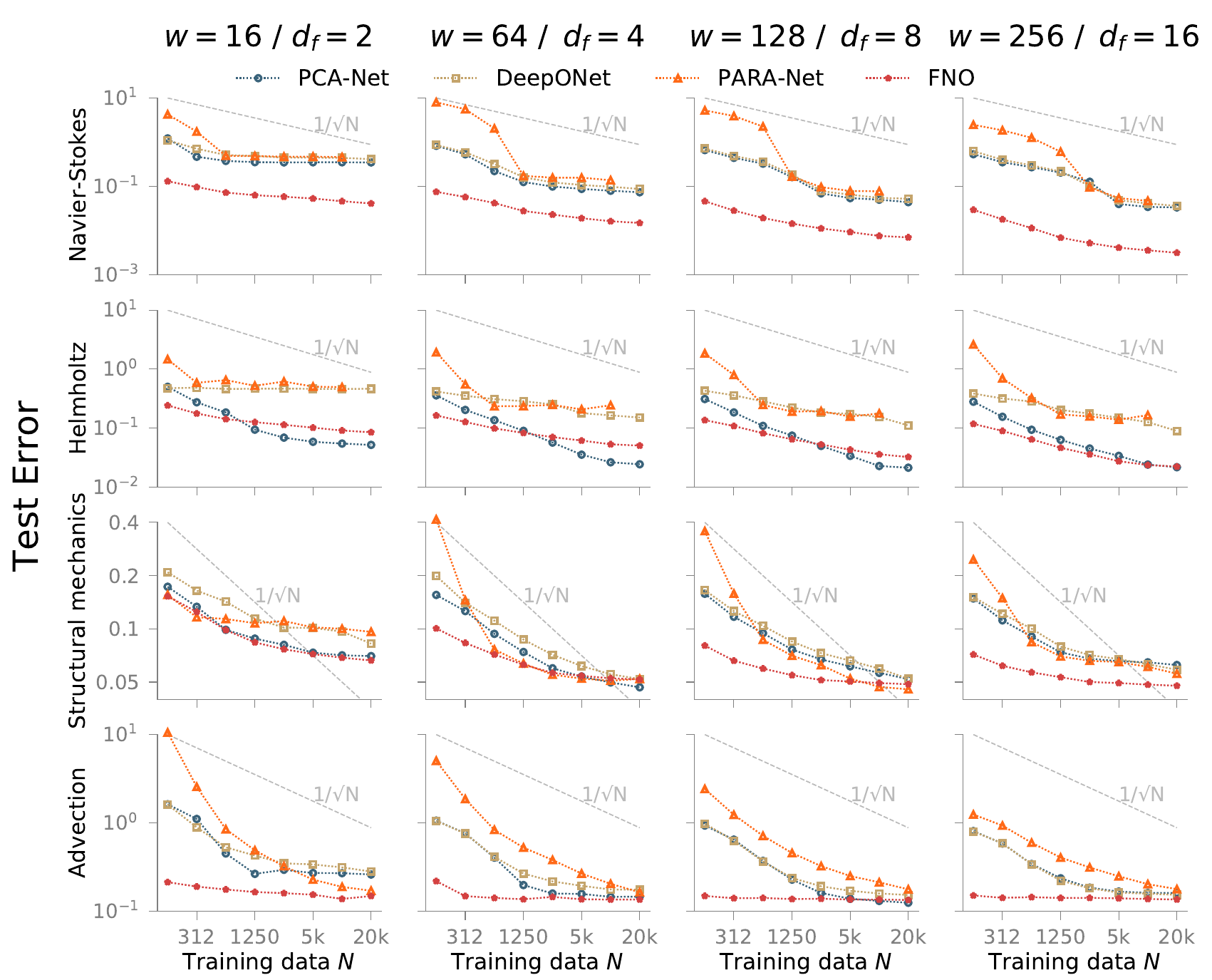}\\
     \caption{Test error vs.\ training data amount $N$ for Navier-Stokes, Helmholtz, structural mechanical, and advection problems (top to bottom). Network size increases left to right.}
     \label{fig:d-e}
\end{figure}

We now begin our discussion of the cost-accuracy trade-off by focusing on studying
the test error as a function of the training cost as measured by the used training data volume.  \Cref{fig:d-e} plots the test error vs.\ training data volume for each of our test problems across different network sizes and architectures, and plots the Monte Carlo rate $\mathcal{O}(N^{-\frac12})$~\cite{caflisch1998monte} for reference.

We contrast the error-training data behavior of the two smooth problems (Navier-Stokes and Helmholtz, with that of the two problems with discontinuous outputs (the structural mechanics and advection problems. 
The paper \cite{de2021convergence} shows that for neural network learning of linear operators the Monte Carlo rate is obtained when the problem is smooth enough, but that for less smooth problems the reduction with respect to $N$ is slower than the Monte Carlo rate. 
This theoretical result is also reflected in our numerical findings (which, with the exception of
the advection problem, all concern nonlinear operators) where the smooth problems exhibit empirical error-data curves with slopes close to the $\mathcal{O}(N^{-\frac12})$ rate for all but the smallest neural networks. For these small neural networks, the expressive power of the networks limits the convergence of error with respect to the training data. In contrast, for our non-smooth problems, the error-data curves generally exhibit slopes that are worse than the $\mathcal{O}(N^{-\frac12})$ rate, and errors level out for all neural network sizes tested, indicating that the expressive power of the networks is more limiting for these non-smooth cases.

Finally, we note that PARA-Net differs from the other network types tested in that it often yields much higher test errors for low training data volumes than the other networks: this indicates that PARA-Net requires greater volumes of training data to yield good predictions. PCA-Net and DeepONet have similar error-training data curves on all problems except the Helmholtz problem, where the DeepONet trunk struggles to capture the high-frequency oscillations of the solution; we will discuss this in~\Cref{ssec:bases}. FNO generally yields the lowest errors for a given training data volume and for the advection test problem actually appears to need even less data than the smallest data volume tested.

\subsubsection{Accuracy vs.\ Network Size/Expressivity}\label{sssec: acc vs size}
\begin{figure}[h]
     \centering
     \includegraphics[width=0.95\textwidth]{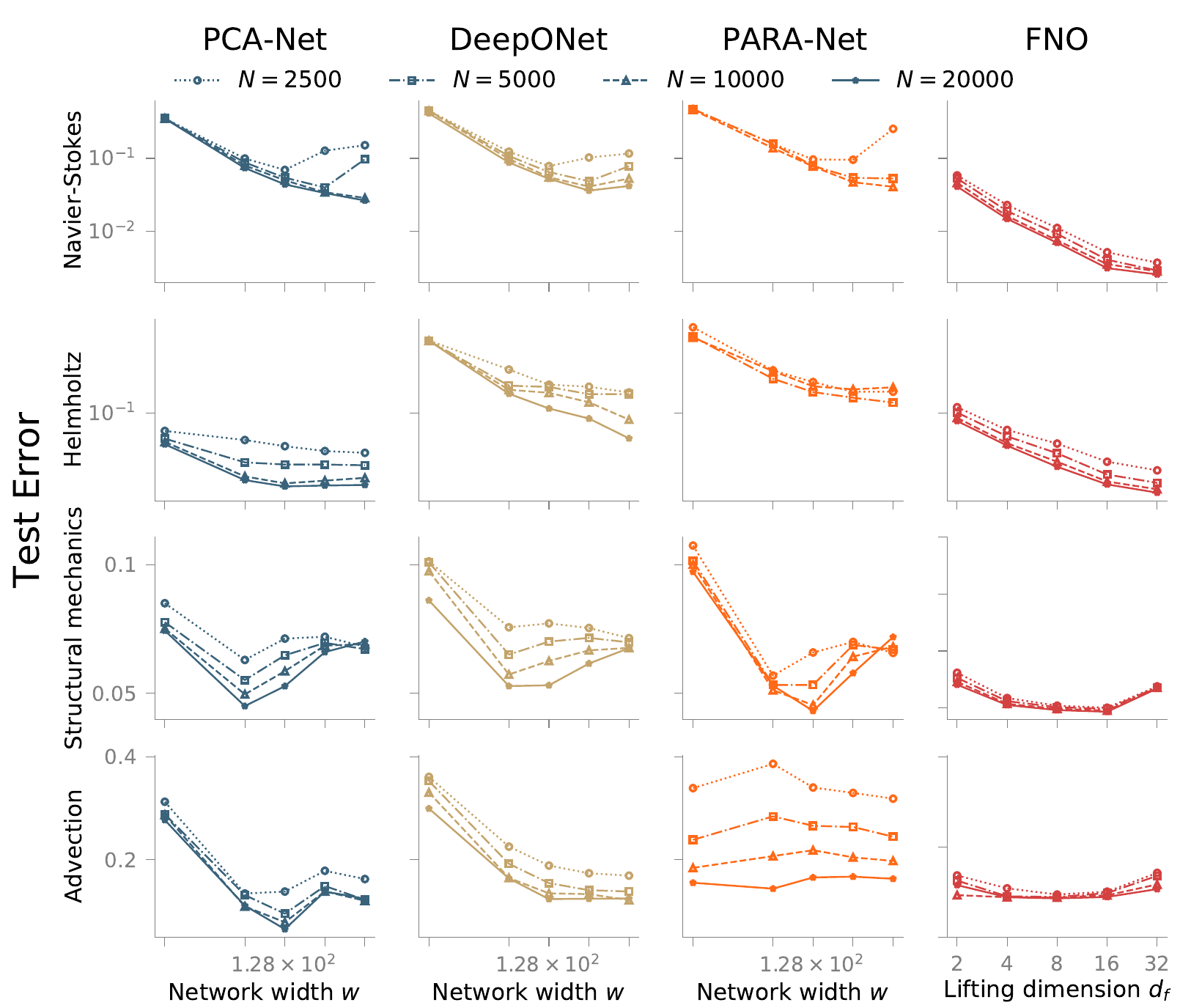}
     \caption{Test error vs.\ network size, as measured by internal layer width $w$ for PCA-Net, DeepONet, and PARA-Net, and as measured by number of channels $d_f$ for FNO. Different lines correspond to different training data volumes, $N$.}
     \label{fig:w-e}
\end{figure}

We next consider the error behavior of our neural network predictions as a function of the network width/number of channels, which is related both to the online evaluation cost of the network and to the network's expressive power. We focus initially here on the trade-off between accuracy (as measured by test error) and expressivity (as measured by number of parameters) and defer discussion of the accuracy-online evaluation cost tradeoff to the next section. 
\Cref{fig:w-e} plots the test error vs.\ network width for each problem and network type.

For the Navier-Stokes test problem, we note that PCA-Net, DeepONet, and PARA-Net all demonstrate overfitting\footnote{The overfitting regime is one where the error grows as the network grows} behavior---that is, error curves that increase with the size of the network---for small training data volumes. However, FNO does not exhibit overfitting behavior for the Navier-Stokes problem, reflecting that FNO's spectral representation of the solution  is particularly suited to this problem. For the Helmholtz test problem, none of the network types exhibit overfitting behavior at any training data volume tested, indicating that $N=2500$ data are sufficient for training the networks for this specific problem.

For our structural mechanics test problem, PCA-Net, DeepONet, and PARA-Net all exhibit overfitting behavior at all training data volumes, whereas FNO exhibits overfitting behavior only at the smaller data volumes. For the advection test problem, both PCA-Net and FNO exhibit slight overfitting behavior at all tested data volumes, whereas DeepONet and PARA-Net appear more robust to this behavior. This may be because DeepONet and PARA-Net both define the output space in terms of a learned neural network, in contrast to PCA-Net and (our specific implementation of) FNO, for which the output space is the result of a linear decomposition or a prescribed set of Fourier bases.

\subsubsection{Accuracy vs.\ Evaluation Cost}\label{sssec: acc v cost}
\begin{figure}[h]
     \centering
     \includegraphics[width=\textwidth]{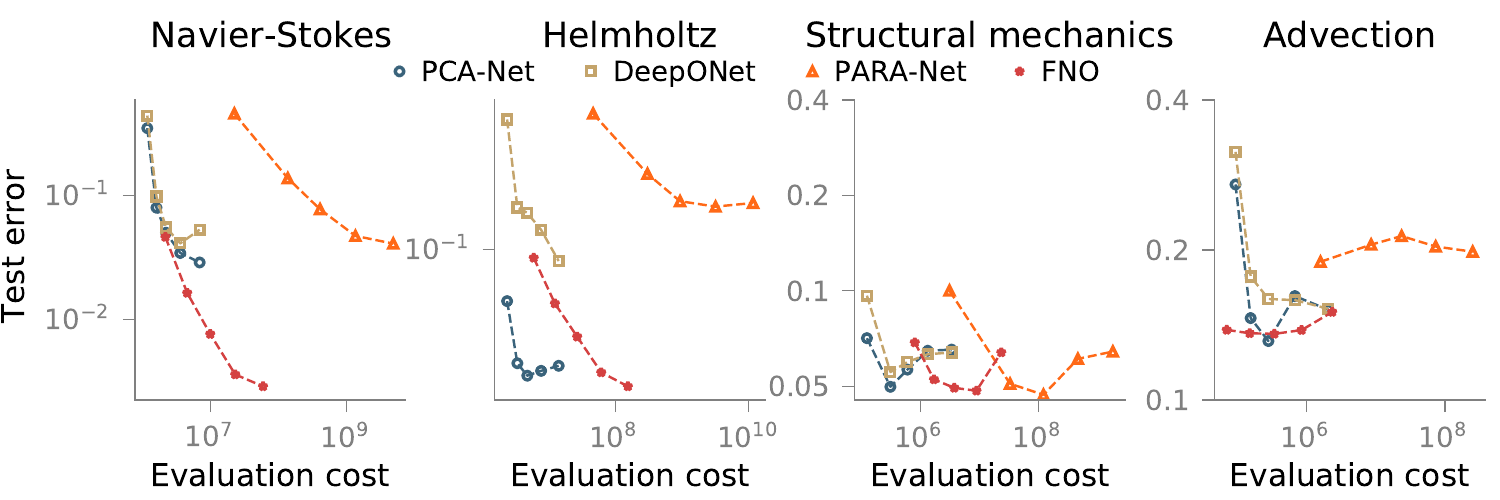}
     \caption{Test error vs.\ evaluation cost for all four test problems for fixed training data volume $N = 10000$. See \Cref{fig:all-error-cost} for error vs.\ evaluation cost curves for all training data volumes tested.}
     \label{fig:c-e}
\end{figure}

Finally, we directly consider the error behavior of our neural network predictions as a function of their online evaluation cost (\Cref{tab: eval complex}). \Cref{fig:c-e} plots the test error vs.\ network evaluation cost for each network type and test problem across the range of network sizes tested, for a fixed training data volume of $N = 10000$. 
\Cref{fig:all-error-cost} contains the complete results for all training data volumes tested; however our main conclusions can be understood from just the $N = 10000$ results.

We note first that for all test problems, the cost of PARA-Net exceeds the cost of other networks by at least an order of magnitude with similar or worse error than the other networks. This is because the network must be evaluated at every spatial point in the discretization of the output function. As such it is not competitive with the other network architectures in the setting where entire output functions are desired; it could however be considered an alternative when outputs at only a few spatial points are desired.

Because our implementation of DeepONet uses a branch network identical to PCA-Net, their evaluation costs are the same for a given network width. Across our test problems PCA-Net generally yields errors similar to those of DeepONet, except for the Helmholtz test problem where the DeepONet trunk struggles to represent the high-frequency features of the solutions that the PCA basis captures. This motivates the use of PCA basis functions as the output space in DeepONet as introduced in~\cite{lu2021comprehensive}.

Finally, our implementation of FNO generally achieves the lowest test errors across all four problems. For the one-dimensional advection equation, the evaluation cost of FNO is similar to that of PCA-Net and DeepONet; however for our two-dimensional examples, the evaluation cost of FNO is greater than that of PCA-Net and DeepONet. This is because the cost of FNO depends on the number of Fourier modes $k_{max}$, which scales exponentially with the spatial dimension. 

Since the cost-accuracy curves for FNO and for PCA-Net/DeepONet occupy different regions of
the cost axis we cannot make definitive conclusions about the relative merits of these three
methods. However, since FNO has the clearest
signal of error decay as a function of cost, we compute
the empirical power law 
$$\textrm{test error}  = a (\textrm{evaluation cost})^{-p}$$
of the FNO. The exponents $p$ are 
$1.129$, $0.7254$, $0.0926$, and $0.002506$ for Navier-Stokes, Helmholtz, Structural mechanics and advection problems, respectively; the differences are presumed to relate to the regularities of the output spaces for these problems.

\eq{
\subsubsection{Out-Of-Distribution (OOD) Generalization}\label{app: out of distribution tests}
Thus far, our numerical studies have exclusively focused on the ability of the neural networks to predict outputs for test inputs drawn from the same distribution as the training data set. In scientific and engineering applications, it is often desirable for a model to be able to accurately predict outputs for inputs outside of the distribution seen in training. We emphasize that the ability of the learned models to generalize to out-of-distribution (OOD) data will be highly dependent on both the problem and the distribution of the unseen data. The overall efficiency of a given surrogate model
will depend on its ability to generalize since this will govern the extent to which
the cost of training can be amortized. Here we provide an initial empirical study of the generalization error to OOD inputs, using the Navier-Stokes and structural mechanics test problems. There is some limited analysis of OOD test error in the linear setting \cite{de2021convergence}, and our numerical experiments in the nonlinear setting add to what is known about this issue.

In both the Navier-Stokes and structural mechanics problems, training input data are drawn from a Gaussian random field with covariance $C$ as described in~\Cref{sssec:NS,sssec:Solid}). To study OOD generalization error, we draw test inputs from a new Gaussian random field with covariance $C' = 4C$, such that the OOD inputs are roughly twice the size of the training inputs.
\Cref{fig:w-d-out-of-d} plots the test error vs.\ training data volume for the largest network sizes tested, with $w = 256$ and $d_f = 16$. As expected, test errors for OOD input data are generally higher than test errors for in-distribution input data. We note that for the smooth Navier-Stokes problem, lower training data volumes lead to higher errors, but the gap between the in-distribution and OOD performance is smaller in these cases. In contrast, for the non-smooth structural mechanics problem, the gap in generalization performance remains large even for smaller data volumes. We finally note that while increasing training data volume improves OOD test error, the convergence rate of the OOD test error with respect to the training data volume is worse than the rate for the in-distribution test error.

\begin{figure}[H]
     \centering
     \includegraphics[width=\textwidth]{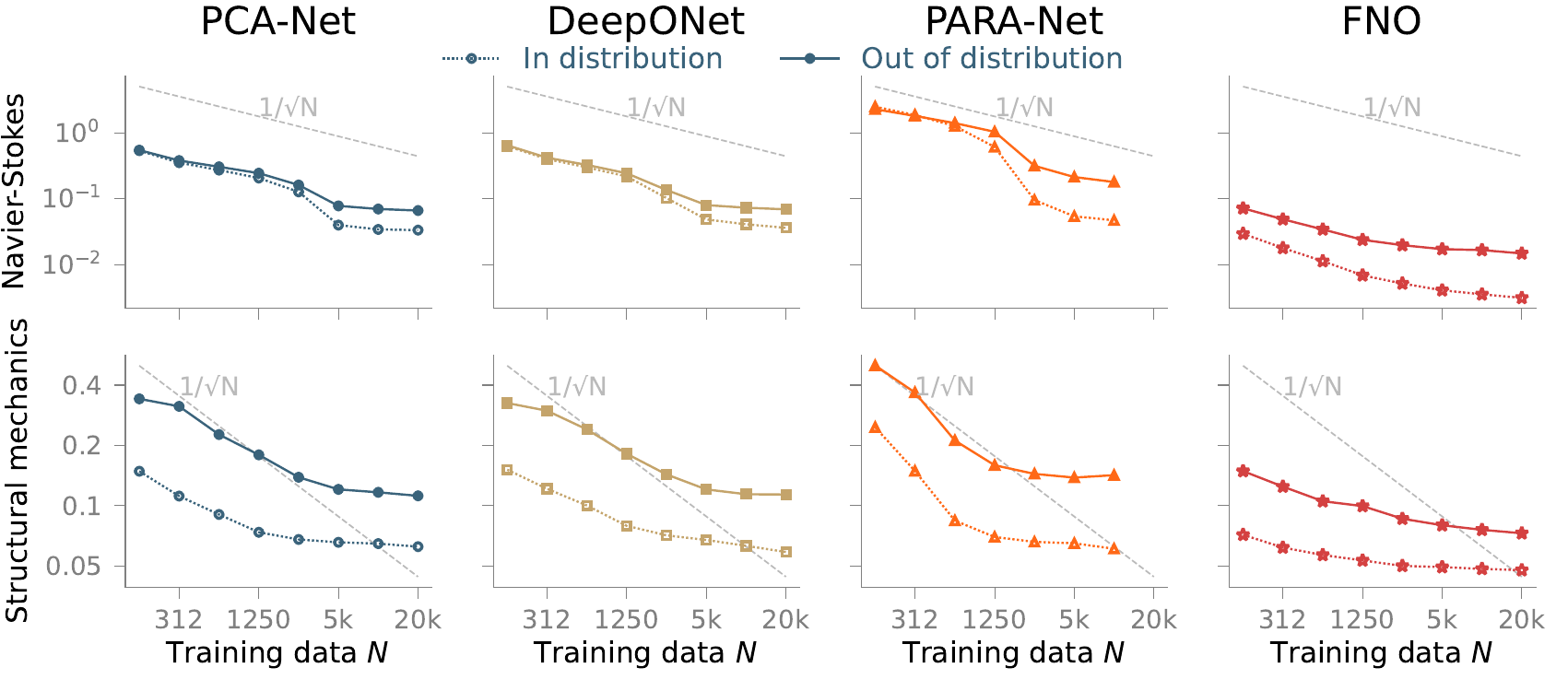}
     \caption{In- vs.\ out-of-distribution test error comparison for Navier-Stokes (top) and structural mechanics (bottom) test problems for largest networks tested ($w = 256$, $d_f = 16$).}
     \label{fig:w-d-out-of-d}
\end{figure}
}

\subsubsection{PCA-Net vs. DeepONet \eq{Output Spaces}}
\label{sssec:output}
\begin{figure}[H]
    \centering
    \includegraphics[width=0.49\textwidth]{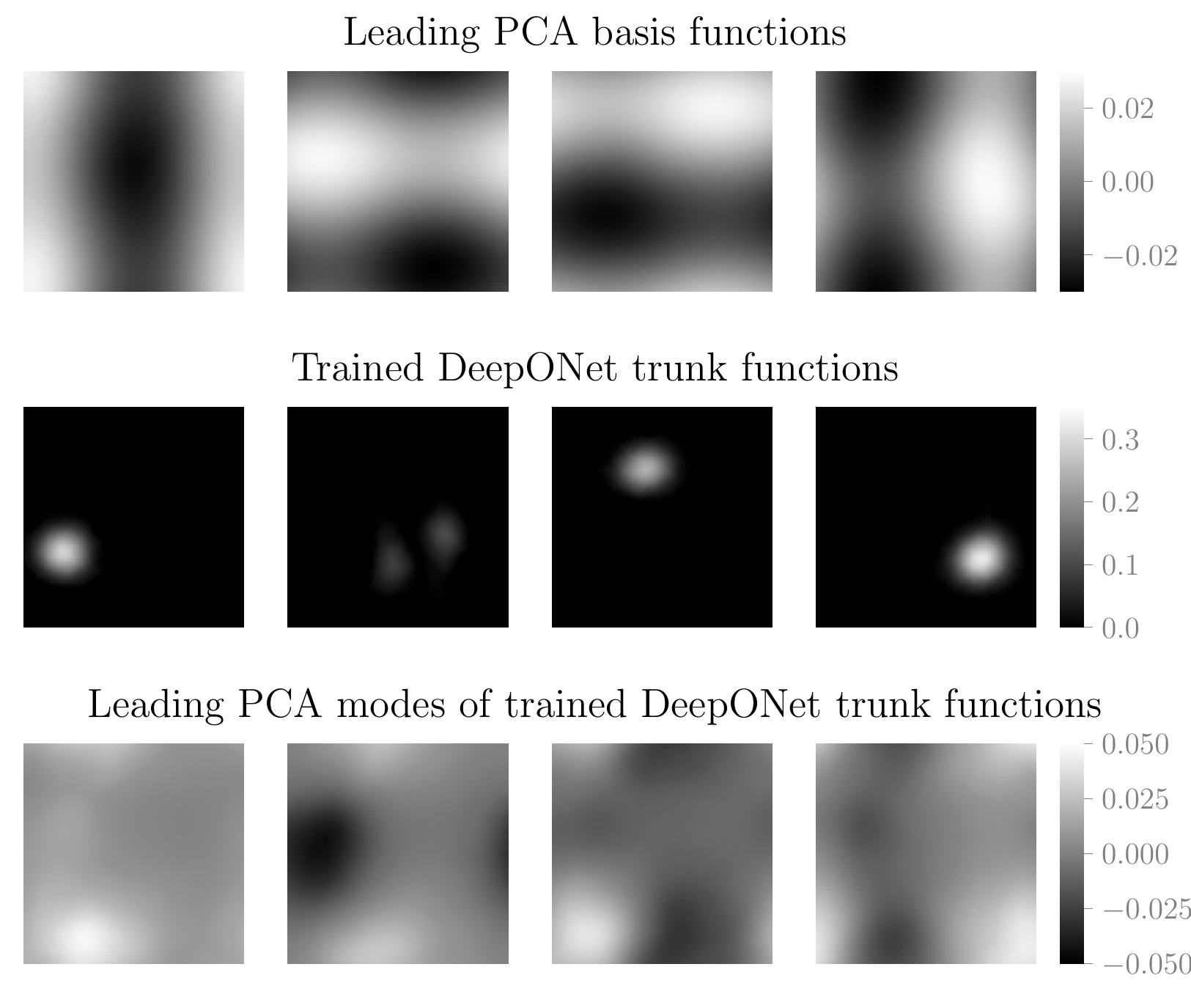}\hfill
     \includegraphics[width=0.49\textwidth]{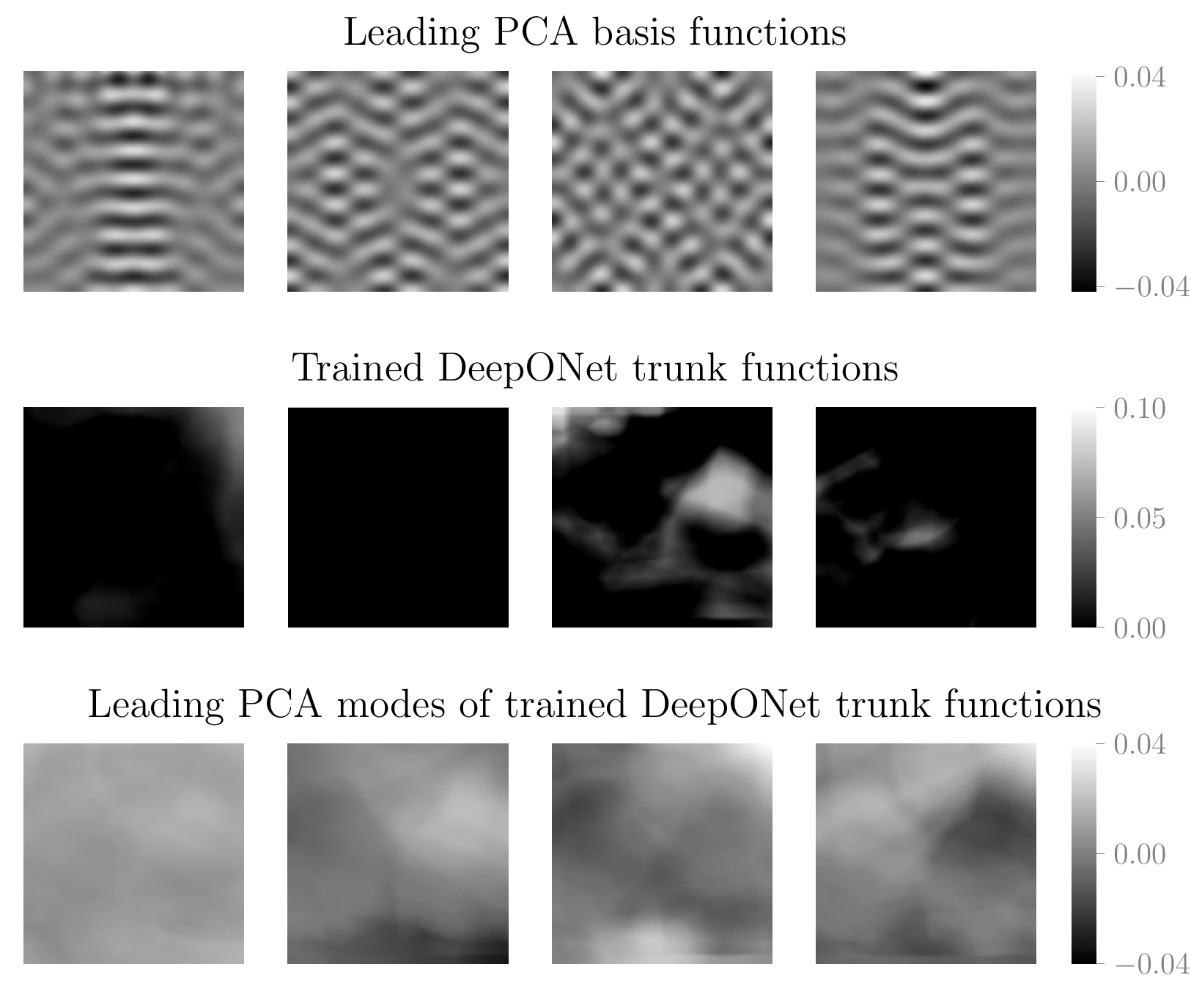}
     \caption{Navier-Stokes (\textbf{left}) and Helmholtz (\textbf{right}) test problems: Comparison of bases used to represent output functions in PCA-Net and DeepONet. Row 1: leading PCA basis functions of output training data. Row 2: Example DeepONet trunk functions after training. Row 3: Leading PCA modes of DeepONet trunk functions after training. }
     \label{fig:smooth-bases}
\end{figure}

\begin{figure}[H]
     \centering
     \includegraphics[width=0.49\textwidth]{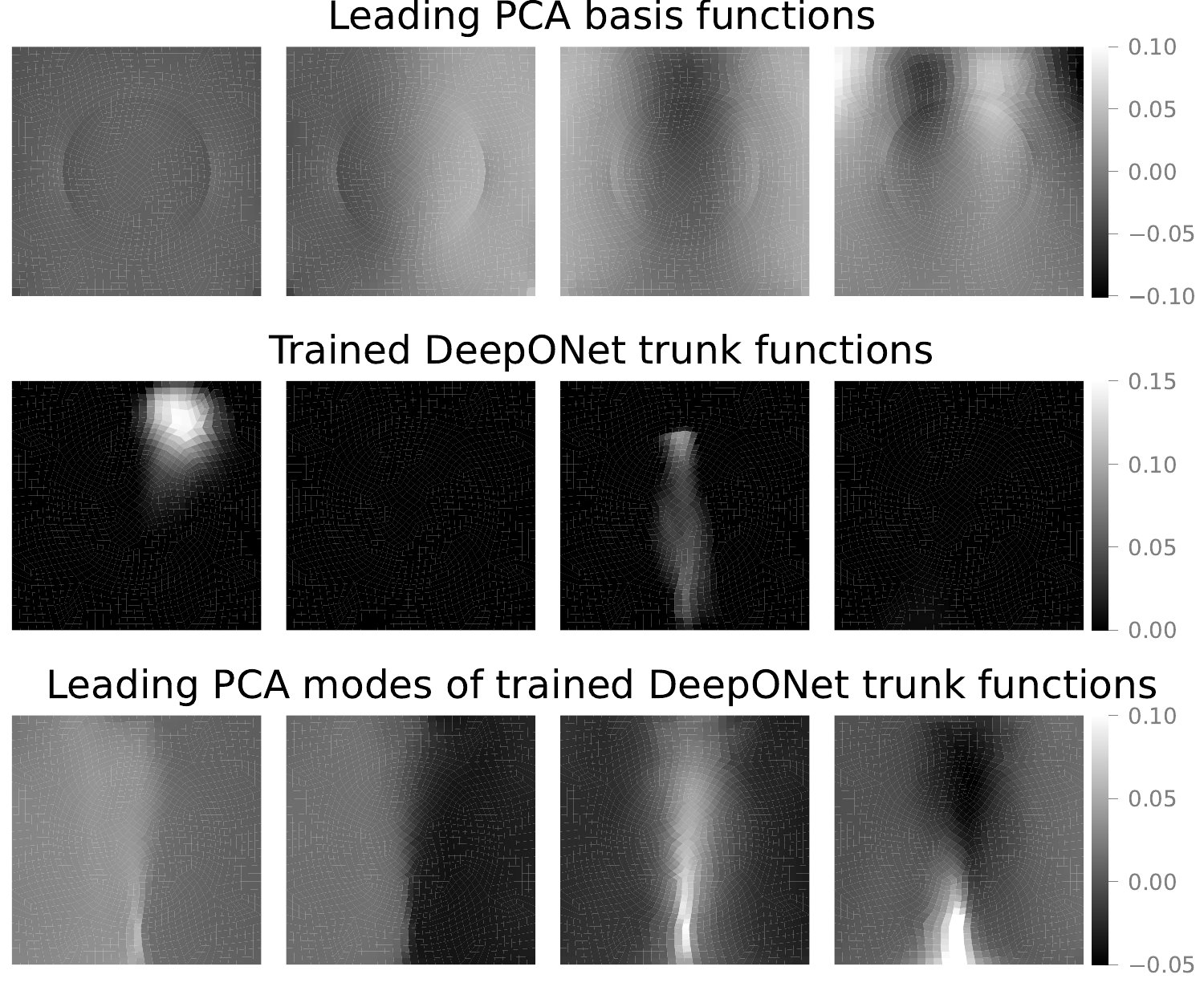}\hfill
     \includegraphics[width=0.49\textwidth]{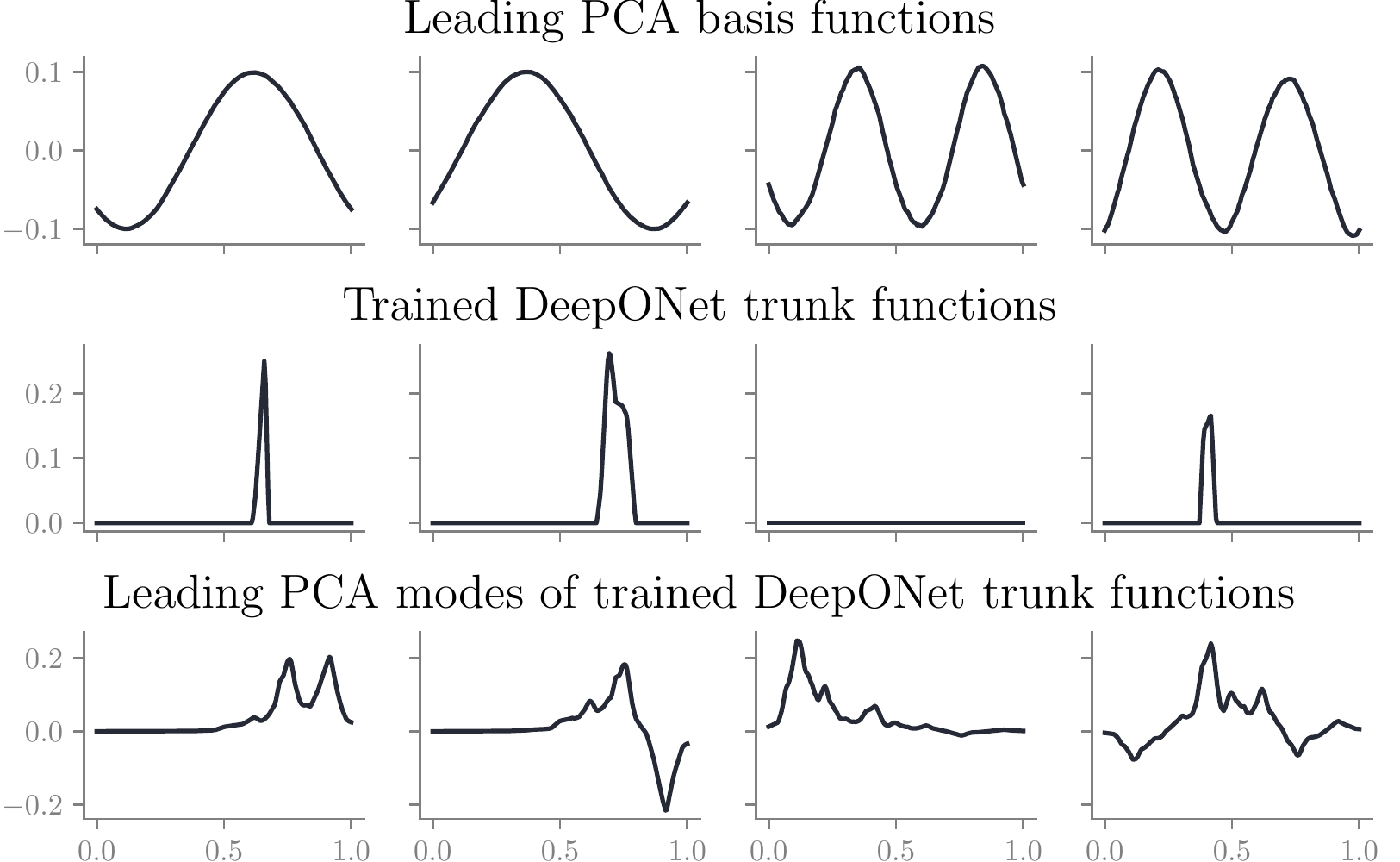}
     \caption{Structural mechanics (\textbf{left}) and advection (\textbf{right}) test problems: Comparison of bases used to represent output functions in PCA-Net and DeepONet. Row 1: leading PCA basis functions of output training data. Row 2: Example DeepONet trunk functions after training. Row 3: Leading PCA modes of DeepONet trunk functions after training.}
     \label{fig:non-smooth-bases}
\end{figure}

Both PCA-Net and DeepONet construct basis functions to represent the output functions, where PCA-Net constructs PCA bases from data directly and DeepONet learns bases from data with the trunk net. It is instructive to compare DeepONet with PCA-Net from the perpsective of the output space. \Cref{fig:smooth-bases,fig:non-smooth-bases} visualize the basis functions used to represent the output functions in PCA-Net and DeepONet for the four test problems we consider. For all test problems, we note that a large number (40\%-70\%) of the DeepONet basis functions after optimization are exactly zero everywhere, despite nonzero random initialization. This may be due to our use of the ReLU activation function, and other activation functions may yield different basis functions. However, our use of ReLU is consistent with the experiments reported in~\cite{lu2021comprehensive}. 
In our visualizations, we overrepresent the non-zero basis functions since these define the output space. Across all examples, we note that the learned DeepONet basis functions after training tend to be local functions, in contrast to the global functions that result when employing in PCA-Net.
 We therefore also compute and visualize the PCA modes of the DeepONet bases.

For the Navier-Stokes problem, PCA-Net and DeepONet achieve similar errors; we attribute this to the fact that the empirical PCA basis of the trained DeepONet trunk functions is similar to that of the true PCA basis of the measure $(\Psid)^\sharp\mu$. However, for the Helmholtz test problem, the PCA basis is able to represent the high-frequency oscillatory nature of the solution, but the empirical PCA basis learned by DeepONet trunk net is not~(\cref{fig:smooth-bases}). 
This explains why PCA-Net achieves lower errors than DeepONet for this problem. This is interesting because in principle if the DeepONet trunk network could learn the PCA basis functions, then DeepONet should be able to achieve similar performance to that of PCA-Net. This raises the question of whether the failure of the trunk to learn the PCA basis functions is due to a lack of expressivity of the trunk or due to optimization error. In our implementation, the DeepONet trunk network only has two hidden layers. However, increasing the number of hidden layers to three actually worsened results in our tests and thus we report the results for two hidden layers. Thus it may be the case that optimization error plays a larger role in the DeepONet error for the Helmholtz problem than for the other three problems.

\label{ssec:bases}

\subsubsection{\eq{Role Of Optimization}}\label{sssec: role of opt}
\begin{figure}[h]
     \centering
     \includegraphics[width=0.9\textwidth]{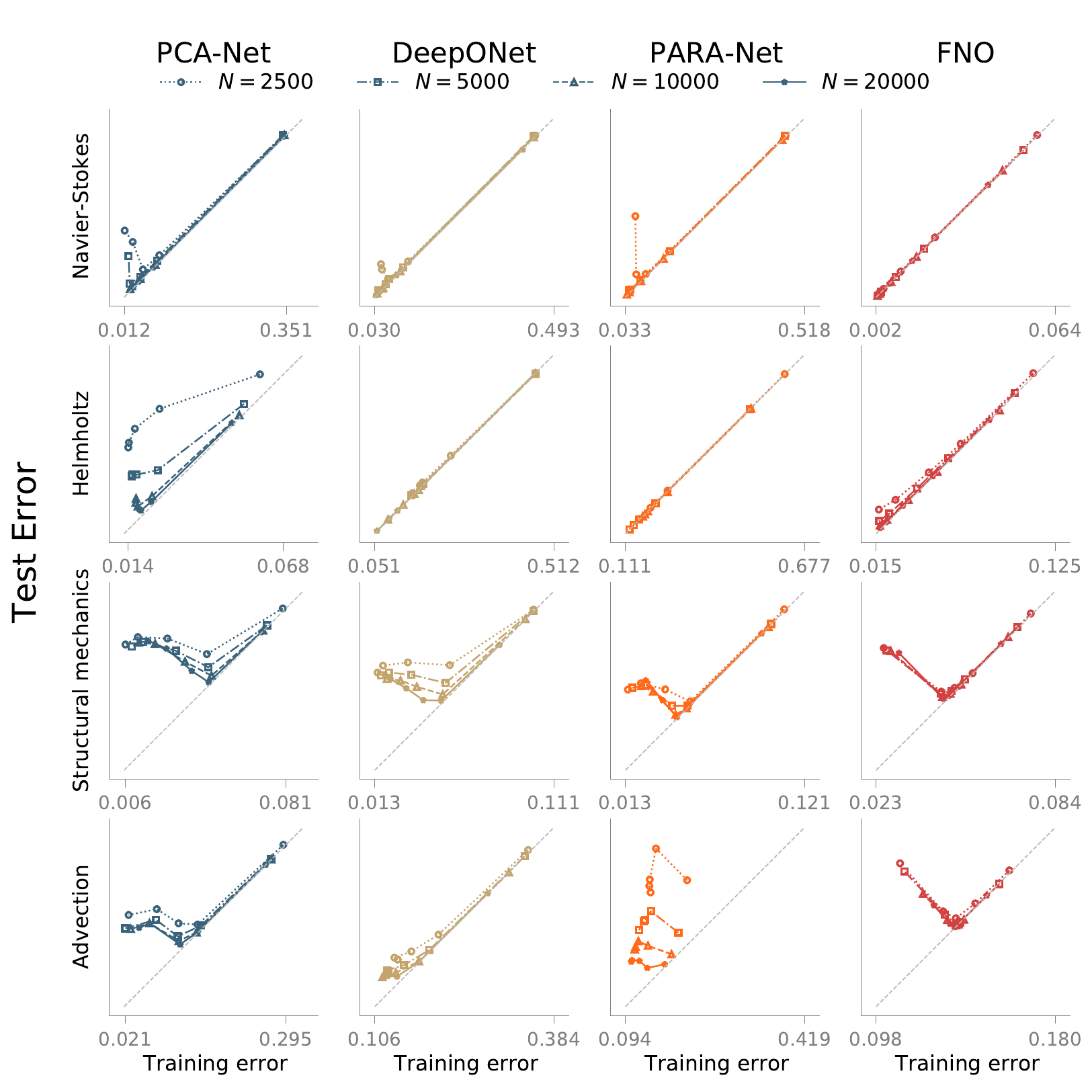}
     \caption{Test error vs.\ training error for all four test problems and all four neural network architectures. Different linestyles correspond to different training data volumes, $N$. The general increase in training error from left to right corresponds with increasing network size as measured by layer width $w$ for PCA-Net, DeepONet, and PARA-Net, and number of channels $d_f$ for FNO. Gray dashed lines have slope 1 and intercept 0.}
     \label{fig:test-v-train}
\end{figure}

\eq{We first offer some comments on the performance of the stochastic gradient optimizer employed in training the neural networks. \Cref{fig:test-v-train} plots the test error vs.\ training error for all four test problems and all four neural network architectures. The different line styles correspond to different training data volumes. The points on each line from left to right correspond to increasing network sizes.
We view the comparison between test and training errors as an indicator of the performance of the optimization. }

\clearpage
\eq{In all cases, test errors are larger than or approximately equal to training errors. Test errors that are noticeably larger than the training error (above the dashed gray lines) are an indication that the optimization has overfit to the training data. This effect is noticeable for the smaller training data volumes and for larger networks, which require more data to train effectively. When test errors are approximately the same as training errors (lying close to the dashed gray lines), the optimization has avoided overfitting to data. Given sufficient training data volume, all four network types avoid overfitting in the smooth test problems (Navier-Stokes and Helmholtz), whereas they generally suffer more from overfitting on the discontinuous test problems (structural mechanics and advection). We note that PARA-Net has particularly bad generalization performance on the advection problem, both compared to its performance on other problems and the performance of other network types on the advection problem. We also note that DeepONet is the most resistant to overfitting on the advection problem, compared to the other network types.}

\eq{Finally we note that the training error is non-zero in all cases, indicating an imperfect fit of the neural networks to the training data. This raises questions of whether different optimization choices would be able to drive the training error lower. Such questions are worthy of study and the subject of future work. We note however that in the cases described above where the neural network has overfit to the training data, better optimization results (in the sense of lower training losses) may not correspond to improved generalization performance. }

\section{Conclusions}
We have presented a numerical study of the performance of four different neural network architectures for modeling operators (function-to-function maps) in problems that involve the solution of a PDE. In particular, we compare test errors for the four networks across a range of network sizes and training data volumes. Our results show that PARA-Net is not a competitive approach in this setting. Our results also suggest that for 1D and 2D problems
\eq{in simple geometries} FNO may be the best approach in terms of the cost-accuracy tradeoff, achieving the lowest errors and low-to-intermediate cost. \eq{However, FNO's cost does scale poorly with spatial dimension and its cost may be prohibitive in the $\ge$3D setting, where DeepONet and PCA-Net may be preferable. It is possible that
investment in efficient software implementations, perhaps taking advantage of parallelization, may make FNO tractable in 3D settings~\cite{takamotopdebench}.} 
\eq{Special treatments of FNO for more complex geometries, a class of problems
not considered here, are introduced in ~\cite{lu2021comprehensive,li2022fourier}; without,
and maybe even with, such special treatments it appears that the additional flexibility of DeepONet gives it some advantages in such settings; a careful complexity study would be of value to establish this.}

Our comparison of DeepONet and PCA-Net fixes the PCA-Net architecture to be the same as that of the DeepONet branch network, so that the key difference between the approaches lies in the functions used to represent the output space. In our comparison, PCA-Net generally yields lower or similar errors to those of DeepONet; investigation into the output functions learned by DeepONet shows that DeepONet does not always succeed at learning a good basis in which to express the output; for the
Helmholtz equation the training process drives the output basis functions to highly localized bases
not representative of the output space. \eq{This result suggests the value of
including PCA basis functions in the output representation of DeepONet for some problems.}

In summary, we provide some insights into the relative merits of different operator surrogates, and mechanisms underlying observed behaviors; whilst we do not reach definitive conclusions about the relative merits of PCA-Net, DeepONet and FNO, and indeed we would expect any such conclusions to be problem dependent, our numerical results highlight the need for, and can guide, future rigorous analyses of the complexity (cost-accuracy trade-off) for these methods. \eq{We have focused on computational studies of error as a function of number of parameters in the network (and hence cost of evaluating the operator surrogate) and as a function of data volume (and hence cost of data acquisition). There is some existing research that relates to these questions, and of course much more is required.
\begin{itemize}
    \item Regarding error as a function of number of parameters, there is profound work     for both DeepONet \cite{lanthaler2021error} and FNO \cite{kovachki2021universal} demonstrating     the possibility of beating the curse of dimensionality on some PDE-based operator learning     problems; in these works this is interpreted as sub-exponential scaling of number of parameters with inverse error. For holomorphic maps arising in parametric elliptic PDEs there is a
    body of literature establishing that certain networks (different from those considered in this paper) can beat the curse of dimensionality \cite{schwab2019deep,herrmann2020deep,opschoor2022exponential}. It is important to appreciate that these results simply concern the \emph{existence} of neural networks with stated properties; they do not address how to find them. Extending the work in these five papers to a broader range of problems, and addressing the issue of finding the desired parameter sets, remain open and challenging problems. The numerical experiments reported in this paper can guide such theory.
    
    \item Regarding error as a function of data volume, the sample complexity of the problem, there is recent theoretical work in the setting of learning linear operators \cite{de2021convergence} and it would be desirable, though challenging, to extend this work to nonlinear operator learning problems. The numerical experiments reports in this paper can also guide such theory.
    
    \item The issue of linking parameterization to data volume, to optimize measures of complexity, is not studied theoretically, to the best of our knowledge. Again the experiments in this paper suggest the need for theoretical guidance on this issue.
    
    \item Questions concerning distribution of parameters through complex network architectures are not studied in any detail, with the exception of studies of depth versus width in neural networks \cite{kutyniok2022theoretical}. It will be important to develop further theory in this area. For example, for the FNO, how should parameter distribution between number of channels and number of Fourier modes be chosen? For
    DeepONet how should parameters be distributed between trunk and branch networks?
    These are hard questions with answers that will be consequential for the choice
    and implementation of operator learning surrogates.
\end{itemize}
}

\vspace{0.1in}
\paragraph{Acknowledgments} MVdH was supported by U.S. Department of Energy, Office of Basic Energy Sciences, Chemical Sciences, Geosciences and Biosciences Division under grant number DE-SC0020345 and the Simons Foundation under the MATH + X program, and the corporate members of the Geo-Mathematical Imaging Group at Rice University. DZH was supported by the generosity of Eric and Wendy Schmidt
by recommendation of the Schmidt Futures program. EQ was supported by Caltech's von Karman postdoctoral instructorship and partially supported by the Simons Foundation Award No. 663281 granted to the Institute of Mathematics of the Polish Academy of Sciences for the years 2021-2023.
AMS was supported by the Office of Naval Research (award N00014-17-1-2079)
and by the Air Force Office of Scientific Research (MURI award number FA9550-20-1-0358 -- Machine Learning and Physics-Based Modeling and Simulation).
The authors thank Mike Kirby, Shibo Li, Zongyi Li,  Lu Lu, Michael Penwarden, Shandian Zhe, and Nicholas Nelsen for helpful comments on an earlier draft.


\appendix

\section{Pointwise Inputs To DeepONet}
\label{sec:pw}

We accommodate the case of pointwise evaluations by writing
$$L_k u = \langle u(x_{\ell}),a_{m} \rangle_{\R^{d_i}}$$
to ensure a collection of real-valued linear functionals on $\cU.$
Note that $k$ is doubly-indexed:
$k=(\ell,m)$ is a multi-index over a set with cardinality defining $d_u;$
the $\{a_{m}\}_{m=1}^{d_i}$ are canonical unit vectors in $\R^{d_i}$, in order to pick-out real-valued functionals,
and the $\{x_{\ell}\}_{\ell=1}^{\ell'}$ denote the locations
of the pointwise evaluations. Thus $d_u=d_i \times \ell'.$
To unify the notation with the PCA input case
we may then (abusing notation) relabel
to index over $k \in \{1,\cdots, d_u\}.$
In this setting it is simplest to think of $\cH=L^2(D_u;\R^{d_i})$
and $\cU=C(D_u;\R^{d_i});$ alternatively $\cU$ may be a RKHS,
such as a Sobolev space of fractional order greater than $d_x/2$,
which is compactly embedded into $C(D_u;\R^{d_i}).$


\section{Complexity Analysis}
\label{app: complexity}

\subsection{Network parameter complexity}\label{app: param complex}
Here we derive the number of parameters of each of the four neural network formulations as a function of the network input and output dimensions and of the network width $w$ (or for FNO, the number of features $d_f$).
Note that a standard fully-connected nonlinear or linear layer from $\R^n\to\R^m$ has $nm$ weights and $m$ biases (the `linear' layer is technically an affine transformation).

PCA-Net consists of an initial nonlinear layer from $\R^{d_u}\to\R^w$, two internal nonlinear layers from $\R^w\to\R^w$, and a final linear layer from $\R^w\to\R^{d_v}$. The parameter complexity for PCA is thus $2w^2 + w(d_u+d_v) + 3w + d_v$.  
DeepONet has two networks, a branch and a trunk. The branch network has the same complexity as PCA-Net. The trunk network consists of an initial nonlinear layer from $\R^{d_y}\to\R^w$, two internal nonlinear layers from $\R^w\to\R^w$, and a final linear layer from $\R^w\to\R^{d_vd_o}$. The parameter complexity for DeepONet is thus $(2w^2 + w(d_u+d_v) + 3w+ d_v)  + (2w^2 + w(d_y + d_vd_o) + 3w + d_vd_o) = 4w^2 + w(d_u+d_v+d_y + d_vd_o) + 6w + d_v + d_vd_o$. 
PARA-Net has an initial nonlinear layer from $\R^{d_u + d_y}\to\R^w$, two internal nonlinear layers from $\R^w \to\R^w$, and a final linear layer from $\R^w\to\R^{d_o}$. The parameter complexity for PARA-Net is therefore $2w^2 + w(d_o+d_u+d_y) + 3w + d_o$.
FNO has an initial lifting layer that lifts the input in $\R^{d_i}$ at each spatial point to channels in $\R^{d_f}$ at each spatial point. In our implementation, the lifting is a linear layer from $\R^{d_i}\to\R^{d_f}$. Similarly, the final projection layer of FNO is a pointwise linear layer from $\R^{d_f}\to\R^{d_o}$. There are three internal Fourier operators from $\R^{N_pd_f}\to\R^{N_pd_f}$. 
In each layer, there is a linear map from $\R^{d_f}\to\R^{d_f}$ that is applied pointwise \eq{with parameter complexity $d_f^2$.}
For each of the $k_{max}$ Fourier modes, there are $d_f^2$ parameters in the linear transformation $P_l$. Each internal layer thus has parameter complexity $d_f^2 k_{max}$. The total parameter complexity for our implementation of FNO is therefore 
 $d_fd_i + d_f + d_fd_o+d_o + 3d_f^2 k_{max} + 3d_f^2$. 

The total number of hyperparameters used in each network for each of the four test problems we consider are tabulated in \Cref{tab: pareq all}.

\begin{table}[h]
\centering
\begin{tabular}{l|c c c c c | c c c c c }
     & \multicolumn{5}{c}{Navier-Stokes} & \multicolumn{5}{c}{Helmholtz} \\
     & \multicolumn{10}{c}{$w$ / $d_f$}\\
     Architecture & 16/2 & 64/4 & 128/8 & 256/16 & 512/32 & 16/2 & 64/4 & 128/8 & 256/16 & 512/32 \\\hline
     PCA-Net &  5680 &41152 & 131456 & 459520 & 1705472&  4816 & 37696 & 124544 & 445696 & 1677824\\
     DeepONet & 7328 & 66176  & 230656  & 854528  & 3281920& 6464 & 62720  & 223744  & 840704  & 3254272\\
     PARA-Net & 5264 & 33344  & 99456  & 329984  & 1184256& 4400 & 29888  & 92544  & 316160  & 1156608\\
          FNO & 1747 & 6973 &27865 &111409 &445537& 1747 & 6973 &27865 &111409 &445537 \\ 
          \\
    & \multicolumn{5}{c}{Structural mechanics} & \multicolumn{5}{c}{Advection} \\
     & \multicolumn{10}{c}{$w$ / $d_f$}\\
     & 16/2 & 64/4 & 128/8 & 256/16 & 512/32 & 16/2 & 64/4 & 128/8 & 256/16 & 512/32 \\\hline
     PCA-Net &  2256 &27456 & 104064 & 404736 & 1595904&  7984 & 50368 & 149888 & 496384 & 1779200\\
     DeepONet & 3904 & 52480  & 203264  & 799744  & 3172352 & 9632 & 75392 & 249088  & 891392  & 3355648\\
     PARA-Net & 1840 & 19648  & 72064 & 275200  & 1074688& 7568 & 42560 & 117888 & 366848 & 1257984\\
     FNO & 1747 & 6973 &27865 &111409 &445537& 163 & 637 & 2521 &10033 &40033
\end{tabular}
\caption{Number of parameters in our implementations of the four neural network architectures for different network size parameters $w$ or $d_f$, for each test problem considered.}\label{tab: pareq all}
\end{table}

\subsection{Network evaluation cost}\label{app: eval complex}
Here we derive the evaluation cost of each of
the four neural networks.
A single standard nonlinear layer $\sigma(Ax + b)$ with $A\in\R^{n,m}$ has cost $2mn + n$, due to the matrix-vector product costing $(2m-1)n$, the vector-vector sum costing $n$, and we approximate the activation function cost as $n$. A single standard linear layer costs $2mn$ since no activation function is applied.
Projecting the solution $u\in\R^{N_p \times d_i}$ on $d_u$ PCA bases has cost $d_u(2N_pd_i-1)$
Recovering the solution $v\in\R^{N_p\times d_o}$ from $d_v$ PCA coefficients has cost $(2d_v-1)N_pd_o$.

PCA-Net has 3 internal layers, the cost is 
 $d_u(2N_pd_i-1) + 2d_uw + 4w^2 + 2d_vw + 3w + (2d_v - 1)N_pd_o.$ 
The DeepONet with the precomputed trunk has the same cost as PCA-Net.
PARA-Net has 3 internal layers. We need to evaluate it at $N_p$ points, and hence its cost is 
$d_u(2N_pd_i-1) + [2(d_u+d_y)w + 4w^2 + 2wd_o + 3w]N_p.$
FNO has an initial lifting layer : $\R^{d_i} \rightarrow R^{d_f}$ with cost $2N_pd_f d_i$, a final projection layer: $\R^{d_f} \rightarrow R^{d_o}$ with cost $2N_pd_od_f$, and 3 Fourier layers.
The cost of the Fourier operator is approximated as $d_f 2\times5N_p  log(N_p) + k_{max}(2d_f^2 - d_f)$, the cost of $\sigma$ is $d_fN_p$, and the cost of $W_l v(x)$ is $N_p(2d_f^2-d_f)$. Hence the cost of FNO is 
$2N_pd_fd_i + 3(10d_fN_p  log(N_p) + k_{max}(2d_f^2 - d_f) + d_fN_p + N_p(2d_f^2-d_f)) +   2N_pd_od_f$.

\section{Error vs.\ cost results}
\begin{figure}[h]
    \centering
    \includegraphics[width=0.8\textwidth]{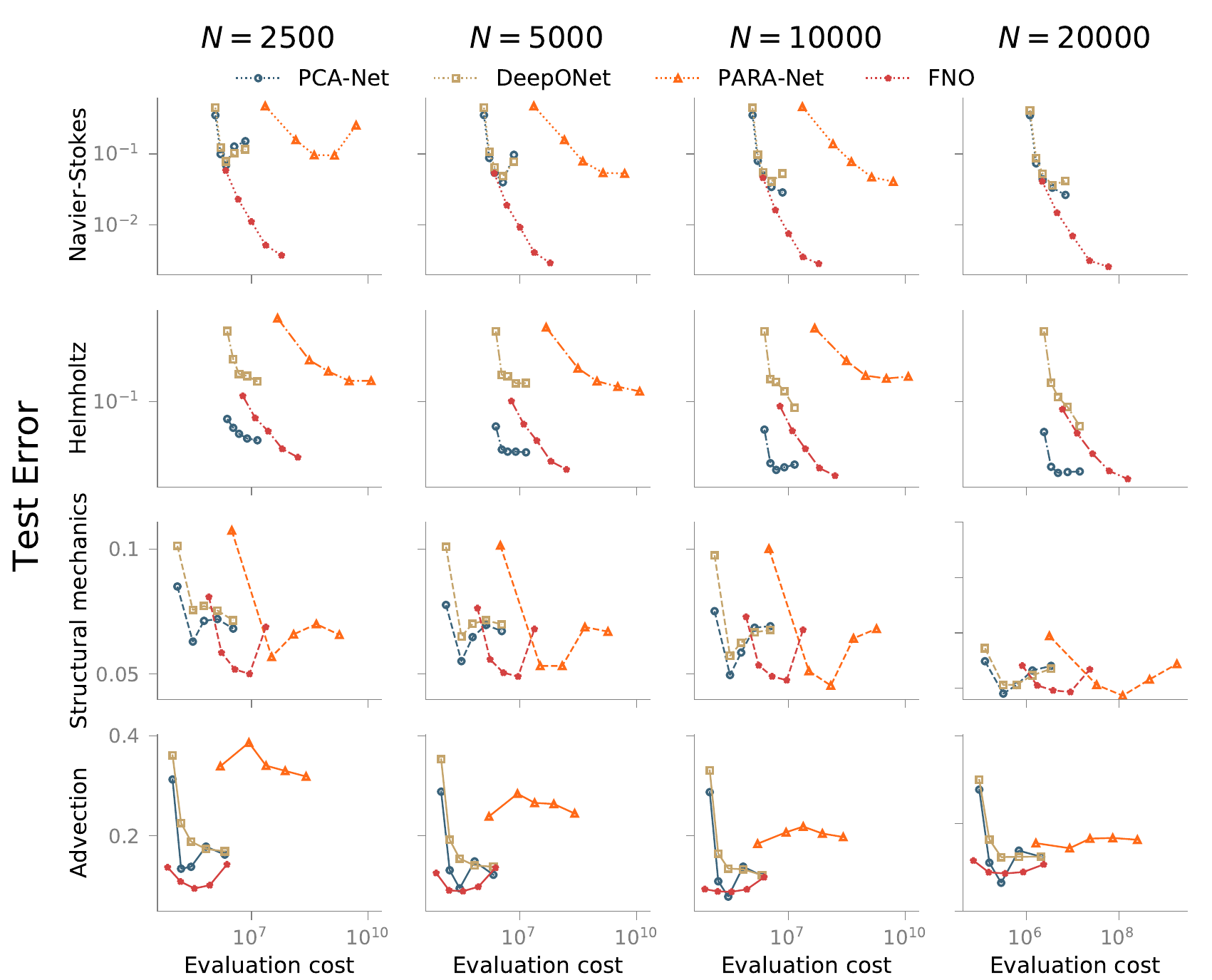}
    \caption{Test error vs.\ evaluation cost for all training data volumes tested.}
    \label{fig:all-error-cost}
\end{figure}

We report in \Cref{fig:all-error-cost} the test error vs.\ online evaluation cost results for all four test problems at all training data volumes tested. The results for each column are similar, and thus only the third column is extracted and reported in \Cref{fig:c-e} in the main text.

\bibliographystyle{unsrt}
\bibliography{references}
\end{document}